\documentclass[11pt]{amsart}
\usepackage{amssymb}
\usepackage{amsthm}
\usepackage{graphicx}
\usepackage{xcolor}
\usepackage[normalem]{ulem}

\theoremstyle{plain}
\newtheorem{Theorem}{{\bfseries Theorem}}[section]
\newtheorem{Lemma}[Theorem]{{ Lemma}}
\newtheorem{Corollary}[Theorem]{{ Corollary}}
\newtheorem{Definition}[Theorem]{{ Definition}}
\newtheorem{Example}[Theorem]{{ Example}}
\newtheorem{Remark}[Theorem]{{ Remark}}

\newcommand{\R}{\mathbb{R}}

\newcommand{\N}{\mathbb{N}}
\newcommand{\Z}{\mathbb{Z}}

\newcommand{\cO}{\mathcal{O}}

\begin{document}
	
	\title[Some Variations of  Transitivity for CR-dynamical systems ]
	{Some Variations of Topological Transitivity for CR-dynamical systems}

\author{Nayan Adhikary}	
	\address{Department of Mathematics, Jadavpur University, Kolkata - 700032, INDIA}
	
	\author{  Anima Nagar}
	\address{Department of Mathematics, Indian Institute of Technology Delhi,
	Hauz Khas, New Delhi 110016, INDIA.}
	
	
	\begin{abstract} We consider the topological dynamics of closed relations(CR) by  studying one of the oldest dynamical property - `transitivity'. We investigate the two kinds of (closed relation) CR-dynamical systems - $(X,G)$ where the relation $G \subseteq X \times X$ is closed and $(X,G, \bullet)$ giving the `suitable dynamics' for a suitable closed relation $G$, where $X$ is assumed to be a compact metric space without isolated points.
		
		$(X,G)$ gives a general approach to study initial value problems for a set of initial conditions, whereas $(X,G, \bullet)$ gives a general approach to study the dynamics of both continuous and quasi-continuous maps.
		
		We observe that the dynamics of closed relations is richer than the dynamics of maps and  find that we have much more versions  of transitivity for these closed relations than what is known for maps. 		
	\end{abstract}

	\maketitle
	
	\noindent{\keywordsname{}:} {\scriptsize{CR-dynamical systems, suitable dynamics, transitivity, weakly mixing, mixing, minimality, strongly transitive, exact transitive, strongly product transitive, locally eventually onto.}}
	
	\tableofcontents
	
	\section{Introduction}
	
	We celebrate thirty years of ``The general topology of dynamical systems'', henceforth cited as \cite{Akin book}, which is the first known work describing the topological dynamics of closed relations on compact metric spaces.
	
	Dynamics is usually represented as differential equations, with trajectories associated to some initial conditions. It is well known that its underlying behaviour be better represented by a difference equation taking into account  the recursive patterns involved, than looking for the associated closed formulas. It then becomes natural to formulate these recursive patterns for a range of initial conditions   among variables in the desired differential equation.  This necessitates the study of dynamics via closed relations blending well with the way we typically think about systems.
	
	After twenty years \cite{Akin book} had an appropriate extension \cite{akin top proc}, which emphasized the use of closed relations to also describe the dynamics of a discontinuous map on a compact metric space.
	
	However, both \cite{Akin book} and \cite{akin top proc} overlook any study of the important and one of the oldest dynamical  property of `transitivity' for dynamics given  by closed relations. Though, in some way such a study has been made recently in \cite{mini, trans, cao}.

	\bigskip
	
Throughout $\mathbb{R}$, $\Z$ and $\mathbb{N}$ stand for the usual space of all real numbers, integers and natural numbers, respectively; and all spaces in the sequel are assumed to be compact metric spaces without isolated points, though our definitions and most results also hold for non-metrizable Hausdorff spaces. 

For any   $A \subseteq X$, we denote its interior by $A^\circ$,  the set of its accumulation or cluster points by $A'$ and its closure by $\overline{A}$. Throughout opene stands for nonempty, open set.

\bigskip

	In this article, we discuss the property of transitivity in its various forms for the dynamics given by closed relations. In section 2, we recall some known concepts and results that will be needed for our later investigation. In section 3, we discuss some known and some new properties of dynamics of closed relations. Finally, in section 4, we study the concept of transitivity for dynamics of closed relations. We define various concepts of transitivity for  CR-dynamical systems,  list some equivalent conditions and study their properties. We find that the list of different versions of transitivity that we get here is more varied than that for maps, and we do not claim that these versions that we study here form an exhaustive list of possibilities.

	\section{Preliminaries}
	
	\subsection{Closed Relations}
	
	We recall the general theory of  closed relations as  described in the epic \cite{Akin book}, supplemented with the exposition \cite{akin top proc}. 
	
	\bigskip
	
	Let $X$ and $Y$ be compact metric spaces. 
	\begin{Definition}\cite{Akin book}
	A \emph{closed relation} $F: X\to Y$ is a closed subset of $X\times Y$ in the product topology such that for every $x\in X$, the set $F(x)=\{y\in Y: (x,y)\in F\}$ is nonempty.	
	\end{Definition}

For every closed set $A\subseteq X,$ $F(A)$ is closed.

For a closed relation $F \subseteq X \times Y$, consider the coordinate projections $\pi_1:X\times Y \to X$ and $\pi_2:X\times Y \to Y$. We note that $\pi_1(F) = X$. We call the relation $F$ \emph{surjective} if $\pi_2(F )= Y$.

\bigskip

 We can regard a continuous function as a special case of closed relation which is identified by the graph of the function. In the case of an arbitrary function $f:X\to Y$, we  write $\overline{f}$ for the closure of the graph of $f$ in the product topology of $X\times Y$. Clearly, $\overline{f}$ is a closed relation and we have $\overline{f}(x)=\{y:(x,y)\in \overline{f}\}$. 
 
 A \emph{selection function} of a closed relation $F$ is a function $f$ such that  for every $x \in X$, $f(x)\in F(x)$. Clearly, $f$ is a selection function of $\overline{f}$. Moreover for the continuous function $f$, we see that $\overline{f}$ is just the graph of $f$ and $f$ is the only selection function of $\overline{f}$. 
 
 \bigskip
 
 As with functions, we can perform composition  of closed relations. Let $F: X\to Y$ and $G:Y\to Z$ be  closed relations. Then $G \circ F=\{(x,z)\in X\times Z: z\in G(y)$ for some $y\in F(x)\}$. A composition of closed relations is closed. Thus when $F$ is a closed relation on $X$, i.e. $X=Y$, we can iterate defining inductively $F^{n+1}=F^n \circ F=F\circ F^n$ for every $n\in \mathbb{N}$. Note that $ y\in F^n(x) $ if and only if there exists $ (n-1) $ points $ x_1,......,x_{n-1} $ such that $ (x,x_1),  (x_i,x_{i+1}), (x_{n-1},y) \in F  $ for every $ i= 1, \ldots, (n-2) $. Also $ F^n(A)= \{ y : y\in F^n(a) $ for some $ a\in A\} $.
 
 The \emph{inverse relation} $F^{-1}: Y \to X$ is defined by $F^{-1}=\{(y,x):(x,y)\in F\}$ and is  closed if and only if $F$ is surjective. For $ B \subseteq Y $ we have   $F^{-1}(B) = \{x \in  X: F(x) \cap B \neq  \emptyset\}$. Again, the composition $ (G \circ F)^{-1} = F^{-1} \circ G^{-1} $. Thus for a closed relation $F^{-1}$ on $X$, we can iteratively define  $F^{-n}=(F^{-1})^n, \ n \in \N$.

 \bigskip
 
We recall some more properties of closed relations. 	If $f:X\to Y$ is a continuous surjection of compact metric spaces then a closed subset $A\subseteq X$ is called \emph{minimal} for $f$ if $f(A)=Y$ and when $B$ is a closed proper subset of $A$ then $f(B)$ is a proper subset of $Y$. A map $f$ is called \emph{irreducible} if $X$ itself is minimal for $f$. This is equivalent to the condition that $f$ is an \emph{almost one-to-one map}, i.e. the set $I_f=\{x\in X: f^{-1}(f(x))=\{x\}\}$ is dense in $X$. A map $f:X\to Y$ is said to be \emph{almost open} if for every $A\subseteq X$ with $A^{\circ}\neq \emptyset$ we have $[f(A)]^{\circ}\neq \emptyset$.

	\begin{Definition}\cite{akin top proc}
		Let  $F:X\to Y$ be a closed relation on compact metric spaces, i.e., $F \subseteq X \times Y$ is closed.  Define $\pi_{1F}:F\to X_1$ and $\pi_{2F}:F\to X_2$ as the restrictions of $\pi_1$ and $\pi_2$ on $F$ respectively.
		
		 Then $F$ is called a \emph{suitable relation} if $\pi_{1F}$ is irreducible and $\pi_{2F}$ is almost open. 
		 
		 Further if $\pi_{2F}$ is also irreducible then $F$ is an \emph{isomorphism} and $F^{-1}$ is also a suitable relation.
		\end{Definition} 
	
	Note that $\pi_{1F}$ is irreducible implies that the suitable relation $F = \bar{f}$, where $f$ is continuous on a dense subset $D \subseteq X$, and $\pi_{2F}$ is almost open implies that for $A \subseteq X$, $A^\circ \neq \emptyset \Rightarrow [F(A)]^\circ \neq \emptyset$.
	
	\begin{Remark}
		For a suitable relation $F \subseteq X\times X,$ and for every opene $ U \subseteq X $ we have $F(U)^\circ \neq \emptyset$ and so $[F^n(U)]^\circ \neq \emptyset$ for all $n \in \N$.	Also for dense $ D \subseteq X $ we have $ F^{-1}(D) $ is also dense in $ X $.
	\end{Remark}
	\begin{Example} \label{sine}
		 We define $F=\{(x,\sin(\frac{1}{x})):x\in [-1,0)\cup (0,1]\}\cup \{(0,y):y\in [-1,1]\} \subseteq [-1,1] \times [-1,1]$. Clearly, $F$ is a closed relation. Suppose that there is a closed set $B$ with $B\subseteq F$ and $\pi_{1F}(B)=[-1,1]$ then $B$ must be equal to $F$. So $\pi_{1F}$ is irreducible. Moreover for every opene set $A$, $[F(A)]^{\circ}$ is nonempty i.e. $\pi_{2F}$ is almost open. Hence $F$ is a suitable relation on $[-1,1]$. 
		 
		 Here $F\circ F=\{(x,sin(\frac{1}{sin(\frac{1}{x})})):x\in [-1,1]\setminus K\}\cup (K\times [-1,1]),$ where $K=\{\frac{1}{k\pi}:k\in \Z\setminus \{0\}\}\cup \{0\}.$ We can similarly  say that $F\circ F$ is suitable.
	\end{Example}
\begin{Example}\label{composition}
	A function $f:[0,1]\to [0,1]$ is defined by $f(x)=x+\frac{1}{2}$ for every $x\in [0,\frac{1}{2}]$ and $f(x)=x-\frac{1}{2}$ for every $x\in (\frac{1}{2},1]$. Take $F=\overline{f}.$ Clearly, $F$ is suitable but $F\circ F=\{(x,x):x\in [0,1]\}\cup \{(0,1),(1,0)\},$ which is not suitable.
\end{Example}

	Hence in general the composition of two suitable relations $F:X \to Y$ and $G:Y \to Z$ need not be suitable. 
	
	\bigskip

	Define  $G \otimes F=\{(x,y,z)\in X \times Y \times Z: (x,y)\in F$ and $(y,z)\in G\}$ and $ONE_G = \{y\in Y : G(y)$ is a singleton set$\}$.
	\begin{Definition}\cite{akin top proc}
	For suitable relations $F:X \to Y$ and $G:Y \to Z$, the \emph{suitable composition} $G\bullet F: X \to Z$ is defined by
	\begin{center}
		$G\bullet F= \pi_{13}[\ \overline{(G \otimes F) \cap (X \times ONE_G \times Z)}\ ]$
	\end{center}
where $\pi_{13}:X \times Y \times Z \to X \times Z$ is the projection to the product of first and third coordinates. 
\end{Definition} 

\begin{Remark}
	Note that if $  F $ and $ G $ are suitable then $ G \bullet F $ is the unique subset of $G \circ F$ that is a suitable relation. In Example \ref{composition}, we  obtain that $F\bullet F=\{(x,x):x\in [0,1]\}$.

\end{Remark}
\bigskip

If $F$ be a suitable relation on $X$, then its suitable iterates are defined by $F^{\bullet 0}=$ identity, $F^{\bullet 1}=F$ and inductively for every $n\in \mathbb{N}$, $F^{\bullet (n+1)}=F\bullet F^{\bullet n}= F^{\bullet n} \bullet F$. When the suitable relation $F$ is an isomorphism, then $F^{-1}$ is also a suitable relation and we can take the suitable iterates   $F^{\bullet {-1}}=F^{-1}$ and inductively for every $n\in \mathbb{N}$, $F^{\bullet -(n+1)}=F^{-1}\bullet F^{\bullet -n} =(F^{\bullet (n+1)})^{-1}$.

Note that if $F$ and $G$ are two suitable relation then $ONE_{G\circ F}$ is a dense subset of $X$ and $G\bullet F$ is the closure of the set $\{(x,G\circ F(x)):x\in ONE_{G\circ F}\}$ (see Theorem 3.14 of \cite{akin top proc}). If $B$ is any minimal relation contained in $G\circ F,$ then it must contain the above set and so from minimality we can say that $B=G\bullet F.$ This proves the uniqueness of $G\bullet F.$

\bigskip

We now consider continuous maps between closed relations.

Let $F:X \to X$, $G:Y \to Y$ and $ H : Y \to X $ be a closed relations between compact metric spaces.  Define 
$$H  \times H : Y  \times Y \to X  \times X \ \text{by} \ \{(y_1, y_2, y_3, y_4) : (y_1, y_3),(y_2, y_4) \in H\}.$$

 The closed relation $ H \times H $ is just the set product of the two relations with the second and third coordinates switched. Note that,
$$(H \times H)(G) = H \circ G \circ H^{-1} \ \text{and} \ (H \times H)^{-1}(F) = H^{-1} \circ  F \circ H.$$

 The continuous surjective map $ h : Y \to X $ maps $ G $ to $ F $ if and only if $ h \circ G = F \circ h $. If $ h $ maps $ G $ to $ F $ then it maps $ G^n $ to $ F^n $ for every  $ n \in \Z $. In this case $h$ is called a \emph{semi-conjugacy}. Furthermore, if $h$ is a homeomorphism and $ h \circ G = F \circ h $ then $h$ is a \emph{conjugacy} between relations $F$ and $G$. Thus if 
 $ h $ and $ \pi_{1G} $ are surjective and  $ \pi_{1F} $ is irreducible then $ (h \times h)(G) = F $.

Moreover, if $F$ and $G$ are suitable relations and $h$ is also almost open then $ h $ maps $ G $ to $ F $ if and only if $ h \bullet G = F \bullet h $, giving a conjugacy between suitable relations. Also,
 $ h $ maps  $G^{\bullet n}$ to $F^{\bullet n}$ for every nonnegative integer
$ n $, and $ h \bullet G^{\bullet n} = F^{\bullet n} \bullet h $.

\bigskip

\subsection{Quasi-continuous maps}

Recently there have been extensive investigations in extending the theory of dynamical systems to maps with some type of discontinuity which is not too far from continuity. In this direction there is a special interest on dynamics of quasi-continuous maps (see \cite{crannell2}). For compact metric spaces  $X$ and $Y$, a function $f:X\to Y$ is called \emph{quasi-continuous} at $x_0\in X$ if for each open neighbourhood $W$ of $f(x_0)$ and each open neighbourhood $U$ of $x_0$, there exists an opene subset $V$ of $U$ such that $f(V)\subseteq W$. If $f$ is quasi-continuous at each point of $X$, then we say that $f$ is quasi-continuous on $X$. Equivalently $f$ is quasi-continuous if and only if for every open set $V$ of $Y$, $f^{-1}(V)$ is quasi-open i.e $f^{-1}(V)\subseteq \overline{(f^{-1}(V))^\circ}$. Let $C_f$ denote  the set of all points of $X$ where $f$ is continuous and likewise $D_f$  denote the set of all discontinuous points of $f$. It is well known that (see \cite{16 crannell}) if $Y$ is second countable and $f$ is quasi-continuous, the set $C_f$ is residual i.e $D_f$ is of first category. Since when $X$ is a compact metric space, then residual sets are precisely those that contains dense $G_{\delta}$ sets, so $C_f$ in such a case contains a dense $G_{\delta}$ set. 

A quasi-continuous almost-open map is called \emph{quopen}.

\begin{Lemma} \cite{akin top proc,crannell}
	For any function $f:X\to Y$ the single-valued set of the closure $\overline{f}$, i.e $ONE_{\overline{f}}$,  is exactly the set of continuity points of $f$. 
\end{Lemma}
The following result presents a nice relation between quasi-continuous maps and closed relation.
\begin{Theorem} \cite{akin top proc} \label{cr-quasi}
	Let $X$ and $Y$ be two compact metric spaces. Then:
	
	$(1)$ If $g:X\to Y$ is a quasi-continuous function and $F$ is the closure of $g$ in $X\times Y$ then $\pi_{1F}$ is irreducible. Furthermore $g$ is continuous at the points of a dense subset of $X$.
	
	$(2)$ If $F\subseteq X\times Y$ is a closed relation with $\pi_{1F}:F \to X$ irreducible and $g:X\to Y$ is a selection function of $F$, then $g$ is quasi-continuous and $F$ is the closure of $g$.
\end{Theorem}

\begin{Theorem}\label{ni}
	Let $F$ be a surjective closed relation on $X$ such that $\pi_{1F}:F \to X$ is irreducible. Then for every opene  $U\subseteq X$ and every $n\in \mathbb{N}$, $(F^{-n}(U))^\circ \neq \emptyset$.
\end{Theorem}
\begin{proof} We prove this result by induction.
	
	Suppose that $g$ is a selection function of $F$. Then from Theorem \ref{cr-quasi}, we can say that $g$ is quasi-continuous and $\overline{g}=F$. Let $U$ be an opene subset of $X$.  Now $F^{-1}(U)\neq \emptyset$ and $\overline{g}=F$, which implies that there is $z\in X$ such that $g(z)\in U$. From quasi-continuity of $g$ we get an opene set $U_1$ such that $g(U_1)\subseteq U$. Hence $U_1\subseteq g^{-1}(U)\subseteq F^{-1}(U)$, which implies that $(F^{-1}(U))^\circ$ is nonempty.
	
	Now let $U_k=(F^{-k}(U))^\circ$ be nonempty for some $k \in \N$. Clearly, $U_k$ is opene and $F^{-1}(U_k)\neq \emptyset$. Then there is $y\in X$ such that $g(y)\in U_k$. It implies that there is an opene set $U_{k+1}$ such that $g(U_{k+1})\subseteq U_k$. Consequently, $U_{k+1}\subseteq g^{-1}(U_k)\subseteq F^{-1}(U_k)\subseteq F^{-1}(F^{-k}(U))$ and so $U_{k+1}\subseteq F^{-(k+1)}(U)$. Therefore $(F^{-(k+1)}(U))^\circ \neq \emptyset$ and so by the induction hypothesis we can conclude that $(F^{-n}(U))^\circ \neq \emptyset$ for all $n\in \mathbb{N}$.      
\end{proof}

\begin{Remark}
	For a suitable  relation $F$ on $X$ and  for every opene $U,$ we have $(F^n(U))^\circ \neq \emptyset$ for all $n\in \mathbb{Z}$. 
\end{Remark}

\bigskip

\subsection{Transitivity for continuous maps}

A common framework in the study of dynamical systems assumes that the phase space is compact metric space and the self mapping is continuous. 
Usually a dynamical system $(X,f)$ is a pair where $X$  a compact metric space,   $f: X \to X$ a continuous self-map with the dynamics given  by iteration. The dynamics of $(X,f)$ is studied via the orbits $\cO(x) = \{x, f(x), f^2(x), \ldots\} $ for all $x\in X$, where $f^n := f \circ f^{n-1}$ for all $n \in \N$ with $f^0$  denoting the identity map. 

\bigskip

	We recall   different versions of topological
transitivity studied in \cite{vt-1}.

$(X,f)$ is called:
\begin{itemize}
	\item[1.] \emph{Topologically Transitive (TT)} if for every opene $U \subseteq X$, $\bigcup \limits_{n=1}^\infty \ f^n(U)$ is
	dense in $X$.
	
	\item[2.]  \emph{Strongly Transitive (ST)} if for every opene $U \subseteq X$, $\bigcup \limits_{n=1}^\infty \ f^n(U) = X$.
	
	\item[3.] \emph{Very Strongly Transitive (VST)} if for every opene $U \subseteq X$ there is a
	$N \in \N,$ such that $\bigcup \limits_{n=1}^N \ f^n(U) = X$.

	\item[4.] \emph{Minimal (M)} if there is no proper, nonempty, closed invariant subset of $X$.

	\item[5.] \emph{Weak Mixing (WM)} if the product system $(X \times X, f \times f)$ is topologically transitive.
	
	\item[6.] \emph{Exact Transitive (ET)}  if for every pair of opene sets $U,V \subseteq X$, \\
	$\bigcup \limits_{n=1}^\infty \ (f^n(U) \cap f^n(V))$ is dense in $X$.
	
	\item[7.] \emph{Strongly Exact Transitive (SET)}  if for every pair of opene sets $U,V \subseteq X$,
	$\bigcup \limits_{n=1}^\infty \ (f^n(U) \cap f^n(V)) = X$.
	
	\item[8.] \emph{Strongly Product Transitive (SPT)} if for every positive integer $k$ the product system $(X^k, f^{[k]})$
	is strongly transitive.

	\item[9.] \emph{Mixing or Topologically Mixing (TM) } if for every pair of opene  sets $U, V \subseteq X$,
	there exists an $N \in \N$ such that $f^n(U) \cap V \neq \emptyset$ for all $n \geq N$.
	
	\item[10.] \emph{Locally Eventually Onto (LEO)} if for every opene $U \subseteq X$,
	there exists $N \in \N$ such that $f^N(U) = X$, and  so $f^n(U) = X$ for all $n \geq N$.
\end{itemize}

It has been  established for dynamics of maps that:

\vskip .3cm

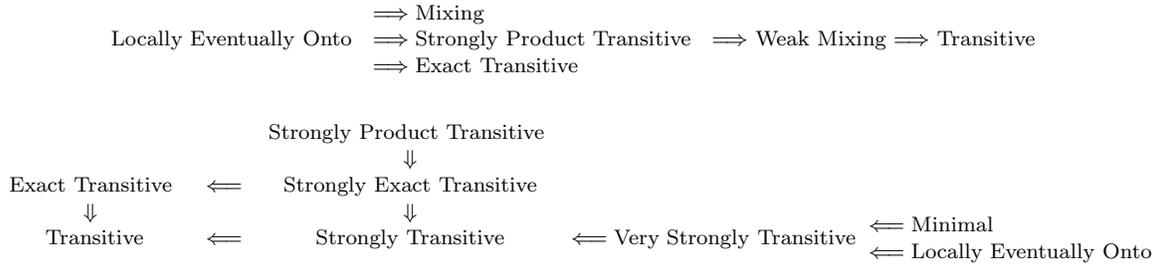
\centerline{\scriptsize{Locally Eventually Onto $\begin{array}{l}
			\Longrightarrow \text{Mixing} \\
			\Longrightarrow \text{Strongly Product Transitive}\\
			\Longrightarrow \text{Exact Transitive}
		\end{array}$
		$\Longrightarrow$ Weak Mixing $\Longrightarrow$ Transitive }}

\vskip .5cm

\centerline{\scriptsize{    $\begin{array}{ccc}
			\text{  } & \text{  } & \text{Strongly Product Transitive  }\\
			\text{  } &  \text{  } & \Downarrow \\
			\text{  Exact Transitive } & \Longleftarrow  & \text{ Strongly Exact Transitive } \\
			\Downarrow & \text{  }  & \Downarrow \\
			\text{ Transitive} & \Longleftarrow  & \text{ Strongly Transitive } \\
			\text{  } & \text{  } & \text{  }\\
			\text{  } & \text{  } & \text{  }\\
			\text{  } & \text{  } & \text{  }\\
			\text{  } & \text{  } &\text{  }
		\end{array}
		\Longleftarrow \text{Very Strongly Transitive} \begin{array}{l}
			\Longleftarrow  \text{Minimal} \\
			\Longleftarrow  \text{Locally Eventually Onto} \end{array}$   }}

Also the reverse implications do  not hold for most of these:

\bigskip

1. Mixing $\nRightarrow$ Strongly Transitive, Exact Transitive or Minimal.

2. Very Strongly Transitive $\nRightarrow$ Minimal.

3. Minimal $\nRightarrow$ Exact Transitive or Weak Mixing.

4. Strongly Product Transitive \& Mixing $\nRightarrow$ Very Strongly  Transitive.

5. Exact Transitive \& Mixing $\nRightarrow$  Strongly  Transitive.

6. Weak Mixing $\nRightarrow$ Mixing.

\bigskip

Recently, Leonel Rito \cite{lrito} showed that Exact Transitive $\nRightarrow$ Mixing.

\bigskip

We include one more form of transitivity viz \emph{point transitive}, which means the existence of a dense orbit, i.e., there exists $x_0 \in X$ with $\overline{\cO(x_0)} = X$. For most systems defined by maps this coincides with the concept of topological transitivity, and so does not merit a seperate study when looking into dynamics of a continuous map.

\bigskip

Let $S$ be a semigroup or monoid acting on a compact metric space $X$. These versions of transitivity were studied for semiflow $(X,S)$ in \cite{vt-2}, where it is shown that some of the implications from above need not hold. Especially it have been shown that the concept of point transitivity is distinct from topological transitivity for dynamics of semiflows.

\bigskip

\section{CR-dynamical systems}

We recall basics from \cite{Akin book,akin top proc} where dynamical systems are presented by closed relations, aided with some more rigourous discussions in the recent articles \cite{mini, trans}. Though, our definitions are  directed to develop our theory on transitivity for CR-dynamical systems.

\bigskip

\subsection{Dynamics of closed relations}

\begin{Definition}
	Let $G$ be a closed relation on a compact metric space $X$, then the pair $(X,G)$ is called a \emph{dynamical system with respect to a closed relation} or briefly a \emph{CR-dynamical system}. 
\end{Definition}

\begin{Definition}
	Let $(X,G)$ be a CR-dynamical system and for every $k\in \mathbb{Z}_+$, the maps $\pi_k:\displaystyle{\prod_{i=0}^{\infty}}X \to X$ denote the $k$-th standard projection of  $\displaystyle{\prod_{i=0}^{\infty}}X$ to its $k^{th}-$ coordinate $X$. Also $G^0$ represents the trivial identity relation. \par
	\begin{enumerate}
		\item For each  $m \in \N$, the set 
		$$\displaystyle{\bigstar_{i=0}^{m}} G =\{(x_0,x_1,x_2,\dots,x_m)\in \displaystyle{\prod_{i=0}^m} X: \ (x_i,x_{i+1})\in G \ \text{for each} \ i\in \{0,1,2,\dots,m-1\} \}$$
		  is called the \emph{$m$-th Mahavier product} of $G$.
		
	\centerline{	\big[We note that the trivial $0$-th Mahavier product is just the identity.\big]}	
		
		 Thus 
		  $$\displaystyle{\bigstar_{i=0}^{m}} G = \displaystyle{\bigcup_{x \in X}} \{(x,x_1,\dots,x_m)\in \displaystyle{\prod_{i=0}^m} X: x_n \in G^n(x), n > 0\}.$$ 
		  
		   \big[Likewise when $G$ is surjective the set $$\displaystyle{\bigstar_{i=0}^{m}} G^{-1} =\{ (x_0,x_1,x_2,x_3,\dots, x_m)\in \displaystyle{\prod_{i=0}^{m}} X: \ (x_{i-1},x_i)\in G^{-1} \ \text{for each} \ i = 1,2, \ldots, m\}$$		   
		    is called the \emph{$m$-th Mahavier product} of $G^{-1}$.\big] 
		  
		  \medskip 
		
		\item The set $$\displaystyle{\bigstar_{i=0}^{\infty}} G =\{ \bar{x} = (x_0,x_1,x_2,x_3,\dots)\in \displaystyle{\prod_{i=0}^{\infty}} X: \ (x_{i-1},x_i)\in G \ \text{ for each} \ i\in \mathbb{N} \}$$ is called the\emph{ infinite Mahavier product} of $G$. 
		
		\medskip
		
	 \big[Likewise when $G$ is surjective the set $$\displaystyle{\bigstar_{i=0}^{\infty}} G^{-1} =\{{\bar{x}}^{-1} = (x_0,x_1,x_2,x_3,\dots)\in \displaystyle{\prod_{i=0}^{\infty}} X: \ (x_{i-1},x_i)\in G^{-1} \ \text{ for each} \ i\in \mathbb{N}\}$$ is called the \emph{infinite Mahavier product} of $G^{-1}$, giving the \emph{inverse limit set} for $(X,G)$.\big] 
	 
	 \medskip 
	 
		\item Let $\bar{x}\in \bigstar_{i=0}^{\infty} G$ and $x_0\in X$. Then $\bar{x}$ is a \emph{forward orbit(trajectory)  of $x_0$} in $(X,G)$ if $\pi_0(\bar{x})=x_0$. 
		
		\medskip
		
		\item The family of all forward orbits(trajectories) of $x_0$ is defined as $T^+_G(x_0)=\{\bar{x}\in \bigstar_{i=0}^{\infty} G : \pi_0(\bar{x})=x_0\}\subseteq \bigstar_{i=0}^{\infty} G$, giving all trajectories starting at $x_0$.  \vspace{0.1cm}
		
	\item A \emph{forward orbit of $x_0$ at $\bar{x}$} is defined as
	\begin{center}
		$ \mathcal{O}^{\oplus}_G(\bar{x})=\{\pi_k(\bar{x}):k \geq 0\}\subseteq X$
	\end{center}
giving the $G$-orbit of $x_0$ associated to $\bar{x}$. \vspace{0.1cm}

\item The \emph{$G$-forward orbit of $x_0$} is defined as
\begin{center}
	$\cO^+_G(x_0) = \bigcup \limits_{n=0}^\infty \ G^n(x_0) = \ 
	\mathcal{U}_G^{\oplus}(x_0)=\bigcup \limits_{\bar{x}\in T^+_G(x_0)} \ \mathcal{O}^{\oplus}_G(\bar{x}) \subseteq X$.
\end{center}

\big[It may be noted here that in general $\mathcal{U}_G^{\oplus}(x_0) \subseteq  \bigcup \limits_{n=0}^\infty \ G^n (x_0)$. But  for an arbitrary relation the reverse containment is not true. 

For example $ G=\{(1,x):x\in [0,1]\} $ is a relation on $ [0,1] $. Here $ x\in G(1) $ for all $x \in [0,1]$, but $\bar{x} \in \displaystyle{\bigstar_{i=0}^{\infty}} G$ with $\pi_0(\bar{x}) = 1$ is a unique point $\bar{x} = \{(1,1,1,\dots) \}$  and so $\mathcal{U}_G^{\oplus}(1) = \{1\}$.  

However for a closed relation $G \subseteq X \times X $ we have $ G(x) $ is nonempty for every $ x\in X $, and so these two sets are equal. \big]

\medskip

\item Let $\bar{x}^{-1}\in \bigstar_{i=0}^{\infty} G^{-1}$ for surjective $G$ and $x_0\in X$. Then $\bar{x}^{-1}$ is a \emph{backward orbit(trajectory)} of $x_0$ in $(X,G)$ if $\pi_0(\bar{x}^{-1})=x_0$. \vspace{0.2cm}

\item The family of all backward orbits(trajectories) of $x_0$ for surjective $G$ is  $T^-_G(x_0)=\{\bar{x}^{-1}\in \bigstar_{i=0}^{\infty} G^{-1} : \pi_0(\bar{x}^{-1})=x_0\}\subseteq \bigstar_{i=0}^{\infty} G^{-1}$,  giving all backward trajectories starting at $x_0$.  \vspace{0.2cm}

\item A \emph{backward orbit of $x_0$ at $\bar{x}^{-1}$} is the set 
		\begin{center}
			$\mathcal{O}^{\ominus}_G(\bar{x}^{-1})=\{\pi_k(\bar{x}^{-1}):k \geq 0\}\subseteq X$.
		\end{center}\par \vspace{0.2cm}

\item For surjective $G$, the \emph{$G$-backward orbit of $x_0$} is defined as
\begin{center}
	$\bigcup \limits_{n=0}^{\infty} \ G^{-n}(x_0) = \cO^-_G(x_0) = \ \mathcal{U}_G^{\ominus}(x_0)=\displaystyle{\bigcup_{\bar{x}^{-1}\in T^-_G(x_0)}} \mathcal{O}^{\ominus}_G(\bar{x}^{-1}) \subseteq X$.
		\end{center}

	\item For a suitable relation $G$, the \emph{suitable $G$-forward orbit of $x_0$} is defined as
	\begin{center}
		$\cO^{\bullet +}_G(x_0) = \bigcup \limits_{n=0}^\infty \ G^{\bullet n}(x_0)  \subseteq X$.
	\end{center}

\item For an isomorphism $G$, the \emph{suitable $G$-backward orbit of $x_0$} is defined as
\begin{center}
	$\cO^{\bullet -}_G(x_0) = \bigcup \limits_{n=0}^{\infty} \ G^{\bullet (-n)}(x_0)  \subseteq X$.
\end{center}

	\end{enumerate}
\end{Definition}

 \begin{Definition}
 	Let $(X,G)$ be a CR-dynamical system and let $A\subseteq X$. 
 	 
(1) The set $A$ is said to be \emph{+invariant} in $(X,G)$ if for each $(x,y)\in G$, $x\in A \implies y\in A$  (It is called $\infty$-invariant in \cite{mini}).
  
 (2)The set $A$ is said to be \emph{-invariant} in $(X,G)$  if for each $(x,y)\in G$, $y\in A \implies x\in A$ (It is called $\infty$-backward invariant in \cite{mini}).
 
 (3) The set $A$ is said to be \emph{invariant} in $(X,G)$  if $A$ is both +invariant and -invariant, i.e. for each $(x,y)\in G, \ x\in A \Longleftrightarrow y\in A$.
 
 \end{Definition}

\begin{Remark}
	Let $(X,G)$ be a CR-dynamical system and $A$ be a +invariant set in $(X,G)$. Let	 $\bar{x}=(x_0,x_1,x_2,x_3, \dots) \in \displaystyle{\bigstar_{i=0}^{\infty}} G$
	 	  and $x_0\in A$. Then $x_k\in A$ for any positive integer $k$. Clearly, $\cO^+_G(x)$  is a +invariant set in $(X,G)$ for every $x\in X$. Similarly $\cO^-_G(x)$  is a -invariant set in $(X,G)$ for every $x\in X$.
\end{Remark}

\begin{Definition}
	Let $(X,G)$ be a CR-dynamical system. Then  $A \subset X$ is called \emph{weakly invariant} in $(X,G)$ if for every $x\in A$, there exists $y\in A$ such that $(x,y)\in G$. (It is called $1$-invariant in \cite{mini}.)
\end{Definition}

\begin{Remark}		
	For every $x\in X$ and every forward orbit $\bar{x}\in T_G^+(x)$, we see that  $\cO_G^{\oplus}(\bar{x})$ is weakly invariant and also $\overline{\mathcal{O}_G^{\oplus}(\bar{x})}$ is a closed weakly invariant subset of $(X,G)$. 	\end{Remark}

If the set $A$ is a closed weakly invariant in $(X,G)$, then we can define a closed relation $G_A=G\cap (A\times A)$ on $A$. Evidently $(A,G_A)$ is a CR-dynamical system.

\begin{Lemma}\label{cl of wi}
	Let $(X,G)$ be a CR-dynamical system and $A$ be a weakly invariant subset of $(X,G)$. Then $\overline{A}$ and $A^{\prime}$ are also weakly invariant in $(X,G).$
\end{Lemma}
\begin{proof}
	Let $A$ be weakly invariant in $(X,G).$ If $A^{\prime}=\emptyset$ then the result is vacously true. Let $x\in A^{\prime},$ then there is a sequence $(x_n)$ in $A$ with $x_n \to x$. Now $G(x_n)\cap A \neq \phi$ for every $n\in \N.$ Let $y_n\in G(x_n)\cap A.$ Then from compactness, we can conclude that $\{y_n\}$ has a cluster point, say $y.$ Consequently, $(x,y)\in G$ and $y\in A^{\prime}$. Hence $A^{\prime}$ is weakly invariant in $(X,G).$
	
	 Similary we can prove that $\overline{A}$ is weakly invariant.	
\end{proof}
\begin{Remark}
	The above Lemma is not true if $X$ is not compact. For example if we consider a closed relation $G$ on $\R$ such that $G=\{(x,\frac{3}{2}):x\leq 0\}\cup \{(x,\frac{1}{x}):x\in (0,1]\}\cup \{(x,\frac{1}{1+x}):x\in [1,\infty)\}.$ The forward orbit $\bar{1}=(1,\frac{1}{2},2,\frac{1}{3},3,\dots)\in \bigstar_{i=0}^{\infty} G.$ It is evident that $0\in \overline{\cO_G^{\oplus}(\bar{1})}$, but $G(0)\cap \overline{\cO_G^{\oplus}(\bar{1})}= \emptyset$.
\end{Remark}

\begin{Definition}
	Let $(X,G)$ be a CR-dynamical system and let $\bar{x}\in \bigstar_{i=0}^{\infty} G$. The set 
	$\omega(\bar{x})= \lim \bar{x} \subseteq X$  is called omega limit set of $\bar{x}$. Note that $\omega(\bar{x}) \neq \emptyset$ for all $\bar{x}\in \bigstar_{i=0}^{\infty} G$.
	
	Again the set 
	$\omega(\bar{x}^{-1})= \lim \bar{x}^{-1} \subseteq X$  is called omega limit set of $\bar{x}^{-1}$. Note that $\omega(\bar{x}^{-1}) \neq \emptyset$ for all $\bar{x}^{-1}\in \bigstar_{i=0}^{\infty} G^{-1}$.
	
	For each $x\in X$, $\omega_G(x)$ will denote the \emph{$G$-omega limit set} of $x$ defined as $\omega_G(x)=\displaystyle{\bigcup_{\bar{x}\in T^+_G(x)}}\omega(\bar{x})$ and \emph{$G^{-1}$-omega limit set} of $x$ denoted as $\omega_{G^{-1}}(x)$ is defined as $\omega_{G^{-1}}(x)=\displaystyle{\bigcup_{\bar{x}^{-1}\in T^-_G(x)}}\omega{(\bar{x}^{-1})}$.
\end{Definition}

\begin{Remark}
	 Note that for each $\bar{x}\in \bigstar_{i=0}^{\infty} G$, we have $\omega(\bar{x})\subseteq  \overline{\mathcal{O}_G^{\oplus}(\bar{x})}$ and $\omega(\bar{x}^{-1})\subseteq \overline{\mathcal{O}_G^{\ominus}(\bar{x}^{-1})}$ when $\bar{x}^{-1} \in \bigstar_{i=0}^{\infty} G^{-1}$.
	 
	 Likewise for every $x \in X$, we have $\omega_G(x)\subseteq \overline{\mathcal{O}^+_G(x)}$ and $\omega_{G^{-1}}(x)\subseteq \overline{\mathcal{O}^-_{G}(x)}$
\end{Remark}

\bigskip

\subsection{Suitable Dynamics of suitable closed relations}
Let $(X,G)$ be a CR-dynamical system, where $G$ is a suitable relation. Then $G^{\bullet n}$ is also suitable for all $n \in \N$ and we can consider orbits given by  suitable composition of $G.$ This defines a new dynamics on $X$ which we call  \emph{suitable dynamics}.\par 
\begin{Definition}
	For a CR-dynamical system $(X,G),$ where $G$ is a suitable relation, if the iterates are defined using  suitable composition, then the resulting dynamical system is called \emph{suitable CR-dynamical system}. We denote it by $(X,G,\bullet)$
\end{Definition}

\begin{Remark}
	The dynamical system defined by a continuous map is a suitable CR-dynamical system.
\end{Remark}

A well known generalization of dynamical systems given by continuous maps is the one given by quasi continuous maps. Let $f$ be a mapping on $X.$ Then the pair $(X,f)$ is called quasi-continuous dynamical system \cite{crannell2} if $f^n$ is quasi-continuous for every $n\in \N.$ In Theorem \ref{qc}, we show that every quasi-continuous dynamical system induces a suitable CR-dynamical system and also present an equivalent condition of the assumption that $f^n$ is quasi-continuous for all $n\in \N.$ 

If $g$ is selection function of the suitable relation $G,$ one can easily observe that $G\bullet G \subseteq \overline{g^2} \subseteq G\circ G.$ The converse implications are not generally true. In Example \ref{composition}, we can see that $F\bullet F=\{(x,x):x\in [0,1]\}$, $\overline{f^2}=\{(x,x):x\in [0,1]\}\cup \{(0,1)\}$ and $F\circ F=\{(x,x):x\in [0,1]\}\cup \{(0,1),(1,0)\}$. We first prove certain lemmas.

\bigskip

Let $f:X\to X$ be any function. Then the set $\mathcal{O}^{\oplus}_f(x)=\{x,f(x),f^2(x),\dots\}$ is the \emph{forward orbit} of $x$ under $f$. Consider the set $C_f^{\infty}=\{x\in X: f^k(x)\in C_f$ for all $k\geq 0\}$. That means if $x\in C_f^{\infty}$, then $f$ is continuous at every point along the orbit of $x$. Moreover, $x\in C_f^{\infty}$ implies that $f^n$ is continuous at $x$ for every $n\in \N.$

\begin{Lemma}\label{2}
	Let $X$ be a compact metric space and $F$ be a suitable closed relation on $X$. Then every selection function $f:X\to X$ of $F$ satisfies the following property:\par  $\blacklozenge$ if $D$ is an open dense subset of $X$ then $(f^{-1}(D))^{\circ}$ is an open dense subset of $X$.
\end{Lemma}

	\begin{proof}
			Let $D$ be an open dense subset of $X$ and $f$ be any selection function of the suitable relation $F$. First we will show that $f^{-1}(D)$ is dense in $X$. Let $U$ be any opene subset of $X$. Then $(F(U))^{\circ}\neq \emptyset$ and so there exists an opene set $V\subseteq F(U)$. Clearly $D\cap V\neq \emptyset$, which implies that there exists $x\in U$ such that $y \in F(x) \cap  D$. Thus $(x,y)\in U\times D$. As $U\times D$ is opene in $X\times X$ and $\overline{f}=F$, there exists $z\in U$ such that $(z,f(z))\in U\times D$. Hence $f^{-1}(D)$ is dense in $X$.\par
			Next from Theorem \ref{cr-quasi}, it is clear that $f$ is quasi-continuous. Now $D$ is opene with $f^{-1}(D)\neq \emptyset$. So $(f^{-1}(D))^{\circ}\neq \emptyset$. Moreover, $f^{-1}(D)$ is quasi-open i.e $f^{-1}(D)\subseteq \overline{(f^{-1}(D))^{\circ}}$. Consequently, $X=\overline{f^{-1}(D)}\subseteq \overline{(f^{-1}(D))^{\circ}}$, which completes the proof.
		\end{proof}

\begin{Corollary}\label{3}
	Let $F$ be a suitable closed relation on a compact metric space $X$. Then for every selection function $f:X\to X$ of $F$, $C_f^{\infty}$ is residual. Hence $C_f^{\infty}$ is dense.
\end{Corollary}
\begin{proof}
	It is clear that $C_f^{\infty}=\displaystyle{\bigcap_{k=0}^{\infty}}f^{-k}(C_f)$. As $f$ is quasi-continuous, $C_f$ is residual. Then applying Lemma \ref{2}, we can conclude that $f^{-k}(C_f)$ is residual for every $k\geq 0$. Hence $C_f^{\infty}$ is residual.
\end{proof}

\begin{Corollary}\label{suit iteration}
	Let $F$ be a suitable closed relation on a compact metric space $X$ and $f$ be any selection function of $f.$ Then for every $n\in \N,$ the suitable iteration $F^{\bullet n}=\overline{\{(x,f^n(x)):x\in C_f^{\infty}\}}.$
\end{Corollary}
\begin{proof}
	The proof follows from Corollary \ref{3} and Theorem 3.14 of \cite{akin top proc}.
\end{proof}

\begin{Lemma}\label{suit int}
	Let $F$ be a suitable relation on $X$ with selection function $f$ and $U,V$ be two opene sets in $X.$ Then for $n\in \N,$ the followings are equivalent:\par
	\begin{enumerate}
		\item $F^{\bullet n}(U)\cap V \neq \emptyset$.
		\item $(U\cap f^{-n}(V))^{\circ}\neq \emptyset.$
		\item $(U\cap F^{-n}(V))^{\circ}\neq \emptyset.$
	\end{enumerate}
\end{Lemma}
\begin{proof}
	
	$(1)\implies (2)$ Let $U,V$ be two opene in $X$ and $F^{\bullet n}(U)\cap V \neq \emptyset.$ So there exist $x\in U$ and $z\in V$ such that $(x,z)\in F^{\bullet n}.$ Then from Corollary \ref{suit iteration}, we get a point $p\in U$ with $p\in C_f^{\infty}$ such that $f^n(p)\in V.$ Consequently, from continuity of $f^n$ at $p,$ there exists an opene set $W\subseteq U$ such that $f^n(W)\subseteq V$ and so $(U\cap f^{-n}(V))^{\circ}\neq \emptyset.$\par
	$(2)\implies (3)$ is obvious.\par
	$(3)\implies(1)$ Suppose that $U$ and $V$ are two opene sets in $X$ with $(U\cap F^{-n}(V))^{\circ}\neq \emptyset.$ Then there is an opene $W$ such that $W\subseteq U\cap F^{-n}(V).$ From Corollary \ref{3}, we get some point $x\in W\cap C_f^{\infty},$ which gives $(n+1)$ points $\{x,f(x),\dots,f^n(x)\}$ such that $x\in U, f^n(x)\in V$ and each $f^i(x)\in C_f$ for $i=0,\dots,n.$ Hence $(x,f^n(x))\in F^{\bullet n}.$
\end{proof}

\begin{Theorem}\label{qc}
	Let $f:X\to X$ be a quopen map and $F=\overline{f}$. Then for every $n\in \N,$ $f^n$ is quasi-continuous if and only if $F^{\bullet n}=\overline{f^n}.$
\end{Theorem}
\begin{proof}
	Since $f$ is quasi-continuous and almost open, $F$ is suitable (Proposition 3.9 of \cite{akin top proc}). For $n=1,$ $F^{\bullet 1}=F=\overline{f}$. 
	Then from Corollary \ref{suit iteration}, it is clear that $F^{\bullet n}\subseteq \overline{f^n}.$ Furthermore, from Theorem  \ref{cr-quasi}, $\overline{f^n}$ is minimal and so $F^{\bullet n}=\overline{f^n}.$ \par
	
	Conversely, suppose that for every $n\in \N,$ $F^{\bullet n}=\overline{f^n}.$ Let $x\in X$ and $U, V$ be two opene sets containing $x$ and $f^n(x)$ respectively. Since $(x,f^n(x))\in F^{\bullet n},$ from Lemma \ref{suit int}, we can conclude that $(U\cap f^{-n}(V))^{\circ}\neq \emptyset.$ Therefore $f^n$ is quasi-continuous.
\end{proof}

\bigskip

In the following example we show that in a quasi-continuous system $(X,f),$ $\overline{f^2}$ may not be equal to $\overline{f}\circ \overline{f}$.

\begin{Example}
	Let $T$ be a tent map defined on $[0,1]$ i.e. $T(x)=2x$ if $x\in [0,\frac{1}{2}]$ and $T(x)=2-2x$ if $x\in [\frac{1}{2},1].$ Let us define a function $f:[-1,1]\to [-1,1]$ by $f(x)=T(x)$ for $x\in [0,1]$ and $f(x)=2x+1$ for $x\in [-1,0).$ One can observe that $f^n$ is quasi-continous for every $n\in \N.$ Take the closed relation $F=\overline{f}$. Clearly, $F$ is a suitable relation. Here $F\bullet F=\overline{f^2}$ and $F\circ F= \overline{f}\circ \overline{f} =  \overline{f^2}\cup \{(1,1)\}.$ Note that $F\circ F$ is not suitable.\par
	Observe that if we consider another selection function $g$ in such a way that $g:[-1,1]\to [-1,1]$ defined by $g(x)=f(x)$ for all $x\in [-1,1]\setminus \{0\}$ and $g(0)=1,$ then $g^2$ is not quasi-continuous at $x=1.$
\end{Example}
From this above example, we can conclude that every selection function of a suitable relation may not give rise to a quasi-continuous system. In the following theorem, we present an equivalent condition in that direction.
\begin{Theorem}
	Let $F$ be suitable relation on $X.$ Then every selection function $f$ of $F$ forms a quasi-continuous system i.e., for every $n,$ $f^n$ is quasi-continuous if and only if $F^n=F^{\bullet n}$ for every $n\in \N.$
\end{Theorem}
\begin{proof}
	Let every selection function $f$ of $F$ form a quasi-continuous system. Suppose to the contrary that there exists some $n\in \N$ such that $F^n\neq F^{\bullet n}.$ Then there exists a point $(x,y)\in F^n$ but $(x,y)\notin F^{\bullet n}.$ Consequently, we have $(n+1)$ points $x=x_0,x_1,x_2,\dots,x_n=y$ such that $(x_{i-1},x_i)\in F$ for each $i=\{1,\dots,n\}$. Now let us consider a selection function $g$ of $F$ such that $g(x_{i-1})=x_i$ for every $i=\{1,\dots,n\}.$ Let $U$ and $V$ be two opene subsets containing $x$ and $y$ respectively. Since $g^n$ is quasi-continuous and $g^n(x)=y,$ there exists an opene set $W\subset U$ such that $g^n(W)\subset V.$ From Corollary \ref{3}, there exists $z\in W\cap C^{\infty}_g$ and so $g^n(z)\in V.$ Thus for every opene $ U \ni x $ and  $ V \ni y  $, we get a point $ z\in C^{\infty}_g $ such that $ (z, g^n(z))\in U \times V $.  Hence $ (x, y)\in  \overline{\{(z, g^n(z):z\in C^{\infty}_g\}} $. So from corollary \ref{suit iteration}, we can conclude that $(x,y)\in F^{\bullet n},$ which gives a contradiction.\par
	The converse part is clear from Theorem \ref{qc}.
\end{proof}

\begin{Theorem}\label{quasi dyn}
	Let $(X,f)$ be a quasi-continuous dynamical system with $f$ being a quopen map and $F = \overline{f}$. Then for any opene $U, V\subseteq X$ and $n\in \N,$ $F^{\bullet n}(U)\cap V \neq \emptyset$ if and only if $f^n(U)\cap V\neq \emptyset.$
\end{Theorem}
\begin{proof}
	One part follows from Lemma \ref{suit int}. Conversely, let $U,V$ be opene in $X$ with $f^n(U)\cap V\neq \emptyset.$ Since $f^n$ is quasi-continuous, there exists opene $W\subseteq U$ such that $f^n(W)\subseteq V.$ Hence $(U\cap f^{-n}(V))^{\circ}\neq \emptyset$ and so from Lemma \ref{suit int}, $F^{\bullet n}(U)\cap V \neq \emptyset$.
\end{proof}

From the above theorem, we can conclude that the dynamics of quasi-continuous map $f$ on $X$ can be realized as the suitable dynamics $(X,\overline{f},\bullet)$.

\bigskip

We see that suitable dynamics $(X,G,\bullet)$ is different from closed relation dynamics of $(X,G).$ Infact,  we shall subsequently see that suitable dynamics gives certain interesting properties for transitivity of continuous or quasi-continuous maps, which generally may not hold for closed relation dynamics.  

\bigskip

	\section{Transitivities for CR-dynamical systems}
	
	Throughout,  $X$ is a  compact metric space and has no isolated points and $G \subset X \times X$ is closed.	Together they give the CR-dynamical system $(X,G)$.

\subsection{Point, Strongly Point, Topological and Suitable Transitivity}
Here we consider  forms of transitivity which are equivalent in  dynamical systems given by continuous maps. But may differ for CR-dynamical systems.

\bigskip

In this section, some of our definitions and results overlap with \cite{trans}. 

\begin{Definition}
	Let $(X,G)$ be a CR-dynamical system. Then: \par 
	1.  $(X,G)$ is called \emph{Topologically Transitive} or \emph{Transitive} if for every pair of opene subsets $U,V \subseteq X$, there exists $n \in \N$
	such that $G^n(U) \cap V \neq \emptyset$. (It is called \emph{+transitive} in  \cite{trans})
	
	Note that this is same as saying that for every pair of opene subsets $U,V \subseteq X$, there exists $n \in \N$ with some $(x_0,x_1,\dots,x_{n}) \in \bigstar_{i=0}^{n} G $  such that $x_0 \in U$ and  $x_{n} \in V$. \par 
	2. $(X,G)$ is called \emph{Point Transitive}  if there exists a point $x \in X$ such that $\cO^+_G(x)$ is dense in $X$. (It is called type \emph{3 dense orbit transitive} or \emph{3-DO-transitive} in \cite{trans}).\par 
	3.  $(X,G)$ is called \emph{Strongly Point Transitive} if there exists $x\in X$ and $\bar{x} \in T^{+}_{G}(x)$ such that $\mathcal{O}^{\oplus}_G(\bar{x})$ is dense in $X$. (It is called  type \emph{2 dense orbit transitive} or \emph{2-DO-transitive} in \cite{trans})\par	
	\end{Definition}

	It is clear that strongly point transitive $\Longrightarrow$ point transitive.

\begin{Theorem}
	Let $(X,G)$ be a topologically transitive  CR-dynamical system. Then $G$ is surjective.
\end{Theorem}
\begin{proof}
	Let $(X,G)$ be a topologically transitive  system. Since $G(X)$ is closed, if $y \in X$ with $y \notin G(X)$ then there exists opene $U\subseteq X$, such that $y \in U$ and  $U\cap G(X) = \emptyset$. But this contradicts transitivity, and hence $G(X)=X$.
\end{proof}

\begin{Corollary}
	When $(X,G)$ is point transitive or strongly point transitive then $G$ is surjective.
\end{Corollary}

\begin{Remark}
	Note that topological transitivity of $(X,G)$ implies that $G^{-1}$ is a closed relation. Thus $G^n(U) \cap V \neq \emptyset \iff $
	
	 there exists $(x_0, x_1, \ldots, x_n) \in \bigstar_{i=0}^{n} G$  with $x_0 \in U$ and $x_n \in V \iff $ 
	 
	 there exists $(x_n, x_{n-1}, \ldots, x_0) \in \bigstar_{i=0}^{n} G^{-1}$ such that $x_n \in V$ and $x_0 \in U$ 
	 
	  $ \iff U \cap G^{-n}(V) \neq \emptyset$
\end{Remark}

\bigskip

For a CR-dynamical system $(X,G)$ and opene subsets  $U,V \subseteq X$ we define \emph{the hitting time set $N(U,V)$} as:

\begin{center}
	$ N(U,V)= \ \{n \in \N:(x_0,x_1,\dots,x_n) \in \bigstar_{i=0}^{n} G \ \text{with} \ x_0 \in U \ \text{and} \ x_n \in V\}$ 
	
	$ = \ \{n \in \N: G^n(U)\cap V\neq \emptyset \} \ = \ \{n \in \N : U\cap G^{-n}(V)\neq \emptyset\}$ 
	
	 $ = \ \{n \in \N:(x_n,x_{n-1},\dots,x_0) \in \bigstar_{i=0}^{n} G^{-1} \ \text{with} \ x_0 \in U \ \text{and} \ x_n \in V\}	.$
\end{center}

\bigskip

Similarly for every $x\in X$ and opene $U \subset X$, we can define 
$$N(x,U)=N(\{x\},U)=\{n\in \mathbb{N}: (\displaystyle{\bigcup_{\bar{x}\in T_G^{+}(x)}}\pi_n(\bar{x})) \cap U \neq \emptyset\}.$$ 

Furthermore, for every forward orbit $\bar{x}\in T_G^+(x)$ of $x$ we can define $$N_{\bar{x}}(x,U)=\{n\in \mathbb{N}:\pi_n(\bar{x})\in U\}.$$ 


We skip the trivial proof of the following:
\begin{Theorem}	Let $(X,G)$ be a surjective CR-dynamical system. Then the following statements are equivalent.\par 
	\begin{enumerate}
		\item  $(X,G)$ is topologically transitive.\par
		\item $N(U,V)  \neq  \emptyset$ for every pair of  opene  $U,V \subseteq X$.
		\item For each opene $U \subseteq X$, $\cO_G^+(U) = \displaystyle{\bigcup_{k=0}^{\infty}}G^k(U)$ dense in $X$.\par 
		\item For each opene $U \subseteq X$, $ \displaystyle{\bigcup_{k=1}^{\infty}}G^k(U)$ dense in $X$.\par 
		\item For each opene $U \subseteq X$, $\cO_G^-(U) = \displaystyle{\bigcup_{k=0}^{\infty}}G^{-k}(U)$ dense in $X$. \par 
		\item For each opene $U \subseteq X$, $ \displaystyle{\bigcup_{k=1}^{\infty}}G^{-k}(U)$ dense in $X$.
	\end{enumerate}
\end{Theorem}

\begin{Theorem} \label{sp=>t} For a CR-dynamical system  $(X,G)$, strongly point transitive implies transitive.
\end{Theorem}
\begin{proof} Let $U$ and $V$ be opene subsets of $X$. Since $(X,G)$ is strongly point transitive, there exists $x\in X$ and $\bar{x}=(x,x_1,x_2\dots)\in T_G^+(x)$ with $\cO_G^{\oplus}(\bar{x})$  dense in $X$. Then there is some $n\in \N$ such that $x_n\in U$. Let  $W=V\setminus \{x_0,x_1,\dots,x_n\}$. Then there exists some $m>n$ with $x_m\in W$. Conequently, $G^{m-n}(U)\cap V \neq \emptyset$.
\end{proof}

\bigskip

The concepts of point transitivity and transitivity are mutually independent in a CR-dynamical systems. We present this fact in the following examples:
\begin{Example}
	We consider $X=[0,1]$ and a closed relation $G$ on $X$ is defined by $G=([0,1]\times \{0\})\cup (\{\frac{1}{2}\}\times [0,1])$. Here $(X,G)$ is point transitive as $\mathcal{O}_G^+(\frac{1}{2})=[0,1]$. But it is not transitive as $G^n((0,\frac{1}{2}))\cap (\frac{1}{2},1)=\emptyset$ for all $n\in \mathbb{N}$. Hence recalling Theorem \ref{sp=>t}, $(X,G)$ is not strongly point transitive.
\end{Example}

\begin{Example}\label{ex2}
	Let $T:[0,1]\to [0,1]$ be the tent map defined as $T(x)=2x$ if $x\in [0,\frac{1}{2}]$ and $T(x)=2-2x$ if $x\in [\frac{1}{2},1]$. Let $B:[0,1]\to [0,1]$ be the Baker transformation defined as $B(x)=2x$ if $x\in [0,\frac{1}{2}]$ and $B(x)=2x-1$ if $x\in (\frac{1}{2},1]$. We see that both $T$ and $B$  are point transitive, and pick transitive points  $y_1$ and $y_2$  of $T$ and $B$ respectively. We define a map $f:[0,2]\to [0,2]$ by $f(x)=T(x)$ if $x\in [0,1]$ and $f(x)=1+B(x-1)$ if $x\in (1,2]$. Now we consider the closed relation $G=\overline{f}\cup \{(1,(1+y_2)),(2,y_1)\}$ on $X = [0,2]$. Observe that  $\cO_G^+(y_1)$ and $\cO_G^+(1+y_2)$ are dense in $[0,1]$  and $[1,2]$ respectively. 
	
	We  claim that $(X,G)$ is transitive. Let $U$ and $V$ be two opene subsets of $X$. We consider following possible cases:\par 
	\textbf{a.} Suppose that $U\cap (0,1)\neq \emptyset$ and $V\cap (0,1)\neq \emptyset$. Since $(f^n(y_1))$ is dense in $(0,1)$, there are $n_1< n_2$ such that $u=f^{n_1}(y_1)\in U$ and $f^{n_2}(y_1)\in V$. Consequently, $f^{n_2-n_1}(u)\in V$.\par 
	\textbf{b.} Suppose that $U\cap (1,2)\neq \emptyset$ and $V\cap (1,2)\neq \emptyset$. This case is analogous to $(1)$ as $(f^{n}(1+y_2))$ dense in $(1,2)$.\par 
	\textbf{c.} Suppose that $U\cap (0,1)\neq \emptyset$ and $V\cap (1,2)\neq \emptyset$. We can find $k,n,m$ such that $u=\frac{k}{2^n}\in U$ and $f^m(u)=1$. Consequently, $(1+y_2)\in G^{m+1}(u)$. Then we get $r> (m+1)$ such that $G^r(u)\cap V \neq \emptyset$.\par 
	\textbf{d.} Suppose that $U\cap (1,2)\neq \emptyset$ and $V\cap (0,1)\neq \emptyset$. We can find $k,n,m$ such that $u=\frac{k}{2^n}\in U$ and $f^m(u)=2$. Consequently, $y_1\in G^{m+1}(u)$. Then we get $r> (m+1)$ such that $G^r(u)\cap V \neq \emptyset$.\par
	
	\bigskip
	
	Now we will show that $(X,G)$ is not point transitive i.e for every $x\in [0,2]$, $\mathcal{O}_G^+(x)$ is not dense in $[0,2]$. Here we consider the following cases:\par 
	\textbf{p.} Let $x\in [0,1]$. Then either $x\neq \frac{k}{2^n}$ for all $k,n\in \mathbb{N}$, which implies that $G^m(x)=T^m(x)\in (0,1)$ for all $m\in \mathbb{N}$, or $x=\frac{k}{2^n}$, for some $k,n\in \mathbb{N}$, which also implies that $\mathcal{O}_G^+(x)$  is not dense in $[0,2]$.\par 
	\textbf{q.} Let $x\in [1,2]$ and $x\neq 1+\frac{k}{2^n}$ for all $k,n\in \mathbb{N}$. Then $G^m(x)=f^m(x)\in (1,2)$ for all $m\in \mathbb{N}$ impling that $\mathcal{O}_G^+(x)$  is not dense in $[0,2]$. \par 
	\textbf{r.} Let $x\in [1,2]$ and $x=1+\frac{k}{2^n}$, for some $k,n\in \mathbb{N}$. Then either $1\in \mathcal{O}_G^+(x)$ which implies that $\mathcal{O}_G^+(x) \cap (0,1)=\emptyset$, or $2\in \mathcal{O}_G^+(x)$ which implies that $\mathcal{O}_G^+(x)$   meets $(1,2)$ in atmost finitely many points i.e. it is not dense in $[0,2]$.\par
	Here for every pair of opene sets $U,V\subseteq X$, the hitting time set $N(U,V)$ is infinite but $(X,G)$ is not point-transitive.
	\end{Example}

\begin{Theorem}
	Let $G$ be a suitable relation on $X$. Then $(X,G)$ is transitive if and only if for every pair of opene sets $U,V \subseteq X$, $N(U,V)$ is infinite.
\end{Theorem}
\begin{proof}
	Let $(X,G)$ be a transitive CR-dynamical system and $U,V$ be two opene subsets of $X$. It is clear that $G(X)=X$. Suppose that $N(U,V)$ is finite. Let $N$ be the maximum number such that $N\in N(U,V)$. From Lemma \ref{ni}, we have $W=(G^{-(N+1)}(V))^\circ\neq \emptyset$. Then there is a non-negative integer $k$ such that $U\cap G^{-k}(W)\neq \emptyset$. Consequently, $U\cap G^{-k}(G^{-(N+1)}(V))\neq \emptyset$, which implies that $U\cap G^{-(k+N+1)}(V)\neq \emptyset$. Hence $(k+N+1)\in N(U,V)$ giving a contradiction.
	
	The converse is vacuous.
\end{proof}

\bigskip

Let $\{(X_n,d_n)\}_{n \in \N} $ be a family of compact metric spaces. Then the cartesian product  $ \prod \limits_{n \in \N} X_n $ is also a compact metric space with the metric  
$$d(\mathbf{x},\mathbf{y})=\Sigma_{n=1}^\infty \dfrac{d_n(x_n,y_n)}{2^n[1+d_n(x_n,y_n)]};$$

for $\mathbf{x}=(x_n), \mathbf{y}=(y_n)\in \prod \limits_{n \in \N} X_n$, 
that induces a topology  equivalent to the product topology on $ \prod \limits_{n \in \N} X_n $.

Let $\{(X_n,F_n):n\in \mathbb{N}\}$ be a family of CR-dynamical systems together with epimorphisms $\{h_i:(X_{i+1},F_{i+1})\to (X_i,F_i):i\in \mathbb{N}\}$. Recall that  an \emph{epimorphism} $h:(X,G)\to (Y,F)$ is a continuous surjection with $h\circ G = F\circ h$ for CR-dynamical systems $(X,G)$ and $(Y,F)$. If such map exists, then $(Y,F)$ is called a \emph{factor} of $(X,G)$ and $(X,G)$ is called a \emph{lift or extension} of $(Y,F).$

 The closed relation $\prod \limits_n F_n$ on $\prod \limits_n X_n$ is defined as:

   $$(\mathbf{x},\mathbf{y})\in \prod_n F_n \Longleftrightarrow (x_n,y_n)\in F_n \ \text{for every} \ n\in \N,$$
   
   for $\mathbf{x}=(x_n), \mathbf{y}=(y_n)\in \prod \limits_n X_n$. Note that since each $F_n$ is a closed and hence a compact subset of $X_n \times X_n$, the cartesian product $\prod \limits_n F_n$ is also a compact subset of $\prod \limits_n [X_n \times  X_n]$. Thus, the CR-dynamical system $(\prod \limits_n X_n, \prod \limits_n F_n)$ is well defined.

 Let us consider a subset $X \subseteq \prod \limits_n X_n$ such that $$\mathbf{x}=(x_n)\in X \Longleftrightarrow x_n=h_n(x_{n+1}).$$
 
 We now  observe the following:\par
\begin{itemize}
	\item $X$ is closed in $\prod_n X_n.$ Consider the sequence $(\mathbf{x}_n) \in X$ such that $\mathbf{x}_n \to \mathbf{x}$ in $X$. Let for each $n\in \N,$ $\mathbf{x}_n=(x_k^n)$ and $\mathbf{x}=(x_k).$ Since for each $n, k\in \N,$ $x_k^n=h_k(x_{k+1}^n)$, from continuity of $h_k$ it is clear that $x_k=h_k(x_{k+1}).$ Hence $\mathbf{x}\in X.$
	\item $X$ is weakly invariant in $(\prod \limits_n X_n, \prod \limits_n F_n)$. If $\mathbf{x}=(x_k)\in X,$ for each $k\in \N,$ $x_k=h_k(x_{k+1}).$ Let us first choose some $y_1\in F_1(x_1).$ Clearly, $y_1\in F_1(x_1)=F_1(h_1(x_2))=h_1(F_2(x_2)).$ Then there exists some $y_2\in F_2(x_2)$ such that $y_1=h_1(y_2).$ Suppose that there exists $y_k\in F_k(x_k)$ with $y_{k-1}=h_{k-1}(y_k)$ for some $k>1$. Then $y_k\in F_k(x_k)=F_k(h_k(x_{k+1}))=h_k(F_{k+1}(x_{k+1})).$ Consequently, there exists some $y_{k+1}\in F_{k+1}(x_{k+1})$ such that $y_k=h_k(y_{k+1}).$ Thus we can construct $\mathbf{y}=(y_k)$ with $y_k\in F_k(x_k)$ and $y_k=h_k(y_{k+1})$ for each $k\in \N$. This implies that $(\mathbf{x},\mathbf{y})\in \prod \limits_n F_n$ and $\mathbf{y}\in X.$\par 
\end{itemize}

Let $F = {\prod \limits_n F_n}{|_X}$, the restriction of $\prod \limits_n F_n$ on $X.$

	Now it is evident that $(X,F)$ is a CR-dynamical system. We call $(X,F),$ \emph{the inverse limit} of the family $\{(X_n,F_n):n\in \mathbb{N}\}$ together with epimorphisms $\{h_i:(X_{i+1},F_{i+1})\to (X_i,F_i):i\in \mathbb{N}\}$.\par 
Recall that a property $P$ of a dynamical system is called a \emph{residual property} if the trivial system (A system is called trivial if the space $X$ is a singleton.) satisfies $P$ and the property $P$ is inherited by factors, irreducible lifts and inverse limits.

 We recall  that the continuous, surjective map 
 $h:X\to Y$ is irreducible if and only if for every opene $U\subseteq X,$ there exists opene $V$ in $Y$ such that $h^{-1}(V)\subseteq U.$

\begin{Theorem}\label{semiconju}
	The property of transitivity, point transitivity and strongly point transitivity of CR-dynamical systems  is closed under semi-conjugacy.
\end{Theorem}
\begin{proof} Let $(X,G)$ and $(Y,F)$ be CR-dynamical systems. And $ h : X \to Y $ be  a semi-conjugacy. 
	
	\begin{enumerate}
		\item Let $(X,G)$ be transitive, and $U, V \subseteq Y$ be opene subsets. Then $h^{-1}(U), h^{-1}(V) $ are opene in $X$ and there exists $n \in \N$ such that $G^n(h^{-1}(U)) \cap h^{-1}(V) \neq \emptyset$. Then $h(G^n(h^{-1}(U))) \cap h(h^{-1}(V)) \neq \emptyset$ implying that $F^n(U) \cap V \neq \emptyset$. Thus $(Y,F)$ is transitive.
		
		\item Let $(X,G)$ be point transitive. Then there exists $x\in X$ such that $\bigcup \limits_{n=0}^{\infty} G^n(x)$ is dense in $X.$ Let $U$ be an opene in $Y.$ Then there exists $n\in N$ such that $G^n(x)\cap h^{-1}(U)\neq \emptyset.$ Consequently, $h(G^n(x))\cap U\neq \emptyset,$ which implies that $F^n(h(x))\cap U \neq \emptyset.$ This proves that $(Y,F)$ is point transitive.
		
		\item Let $(X,G)$ be strongly point transitive. Then there exist $x\in X$ and $\bar{x}\in T_G^+(x)$ such that $\cO_G^{\oplus}(\bar{x})$ is dense in $X.$ Let $\bar{x}=(x=x_1,x_2,\dots,x_n,\dots)$. Observe that for each $i,$ $h(x_{i+1})\in h(G(x_i))\subseteq F(h(x_i))$. Clearly, $\overline{h(x)}=(h(x_1),h(x_2),\dots,h(x_n),\dots)\in T_F^+(h(x))$ and then one can easily show that $\cO_F^{\oplus}(\overline{h(x)})$ is dense in $Y.$ Hence $(Y,F)$ is strongly transitive.
	\end{enumerate}
\end{proof}

\begin{Lemma}\label{inverselimit}
	Let $\{(X_n,F_n):n\in \mathbb{N}\}$ be a family of CR-dynamical systems together with epimorphisms $\{h_i:(X_{i+1},F_{i+1})\to (X_i,F_i):i\in \mathbb{N}\}$ and $(X,F)$ be the associated inverse limit. Suppose that $p_k:(X,F)\to (X_k,F_k)$ is the kth coordinate projection. Then each $p_k$ is an epimorphism and for every opene $U \subseteq X$ there exist some $k \in \N$ and an opene $V$ in $X_k$ with $p_k^{-1}(V)\subseteq U.$
\end{Lemma}
\begin{proof}
Obviously each $p_k$ is a continuous surjection and $p_k \circ F \subseteq F_k\circ p_k.$ To prove the converse let $\mathbf{x}=(x_n)\in X$ and $(\mathbf{x},y)\in F_k \circ p_k.$ Take $y_k=y$ and $y_{k-1}=h_{k-1}(y).$ Now $y_{k-1}\in h_{k-1}(F_k(x_k))=F_{k-1}(h_{k-1}(x_k))=F_{k-1}(x_{k-1})$. In this way we get $y_{k-1},\dots,y_1$ with $y_i=h_i(y_{i+1})$ and $y_i\in F_i(x_i)$ for $i=1,\dots (k-1).$ On the other hand $y_k\in F_k(x_k)=F_k(h_k(x_{k+1}))=h_k(F_{k+1}(x_{k+1})).$ This implies that there exists some $y_{k+1}\in F_{k+1}(x_{k+1})$ with $y_k=h_k(y_{k+1}).$ In this way we can choose $y_{k+1},y_{k+2},\dots$ such that for each $i\geq k,$ $y_i\in F_i(x_i)$ and $y_i=h_i(y_{i+1}).$ Evidently $\mathbf{y}=(y_n)\in X$ with $y_k=y$ and also $(\mathbf{x},\mathbf{y})\in F.$ Hence $(\mathbf{x},y)\in p_k \circ F.$ This proves that $p_k$ is an epimorphism.\par 
Let $U$ be an opene in $X.$	Then $U=(\prod \limits_n U_n)\cap X,$ where each $U_n$ is opene in $X_n$ and there exists some $k\in \N$ such that $U_n=X_n$ for all $n>k.$ Take $U_1=W_1.$ From continuity of $h_1,$ we can conclude that $W_2=h_1^{-1}(W_1) \cap U_2$ is opene in $X_2.$ In this process we get opene $W_{i+1}=h_i^{-1}(W_i)\cap U_{i+1}$ in $X_{i+1}$ for $i=1,\dots,(k-1).$ Let $\mathbf{x}=(x_n)\in p_k^{-1}(W_k).$ Since $x_i=h_i(x_{i+1}),$ $x_i\in W_i$ for each $i=1\dots,k.$ This implies that $\mathbf{x}\in U.$ Hence $p_k^{-1}(W_k)\subseteq U.$ 
\end{proof}

\bigskip 

\begin{Theorem}
	Topological Transitivity of closed relations is a residual property.
\end{Theorem}
\begin{proof}
From Theorem \ref{semiconju}, it is clear that topological transitivity is inherited by factors. 

Let $(X,G)$ and $(Y,F)$ be two CR-dynamical systems and $h:(X,G)\to (Y,F)$ be an irreducible epimorphism. Suppose that $(Y,F)$ is topologically transitive and $U,V$ are  opene in $X.$ Then there exist opene $U_1$ and $V_1$ in $Y$ such that $h^{-1}(U_1)\subseteq U$ and $h^{-1}(V_1)\subseteq V.$ Now there is some $n\in \N$ and  $F^{n}(U_1)\cap V_1\neq \emptyset.$ This implies $F^{n}(h(h^{-1}(U_1)))\cap V_1\neq \emptyset.$ Then $h(G^n(h^{-1}(U_1)))\cap V_1\neq \emptyset.$ Consequently, $G^n(h^{-1}(U_1)))\cap h^{-1}(V_1)\neq \emptyset.$ Hence $G^n(U)\cap V\neq \emptyset,$ and $(X,G)$ is transitive.\par 

Let $\{(X_n,F_n):n\in \mathbb{N}\}$ be a family of CR-dynamical systems together with epimorphisms $\{h_i:(X_{i+1},F_{i+1})\to (X_i,F_i):i\in \mathbb{N}\}$ and $(X,F)$ be the associated inverse limit and $p_k:(X,F)\to (X_k,F_k)$ be the kth coordinate projection. Assume that each $(X_i,F_i)$ be topologically transitive. Suppose that $U$ and $V$ are  opene in $X.$ Then from Lemma \ref{inverselimit}, there exist some $k$ and  opene $U_1,V_1$ in $X_k$ such that $p_k^{-1}(U_1)\subseteq U$ and $p_k^{-1}(V_1)\subseteq V.$ Then proceeding as above we can conclude that $(X,F)$ is topologically transitive.
\end{proof}

\begin{Theorem} \label{irrex}
	Point transitivity and strongly point transitivity of closed relations are inherited by  irreducible extensions.
\end{Theorem}
\begin{proof}
		Let $(X,G)$ and $(Y,F)$ be two CR-dynamical system and $h:(X,G)\to (Y,F)$ be an irreducible epimorphism.
	\begin{enumerate}
		\item Suppose that $(Y,F)$ is point transitive. Then there exists $y\in Y$ with $h(x)=y$ for some $x\in X$ such that $\displaystyle{\bigcup_{n=0}^{\infty}} F^n(y)$ is dense in $Y$. Let $U$ be opene in $X.$ Then  there exists an opene $V$ in $Y$ such that $h^{-1}(V)\subseteq U.$ Now there exists $n\in \N$ such that $F^n(y)\cap V= F^n(h(x))\cap V \neq \emptyset,$ which implies that $h(G^n(x))\cap V\neq \emptyset.$ Consequently, $G^n(x)\cap h^{-1}(V)\neq \emptyset.$ Hence $\displaystyle{\bigcup_{n=0}^{\infty}}G^n(x)$ is dense in $X.$ 
		
		\item Suppose that $(Y,F)$ is strongly point transitive. Then there exists a forward orbit $\bar{y}=(y=y_1,y_2,\dots,y_n,\dots)$ such that $\cO_F^{\oplus}(\bar{y})$ is dense in $Y.$ Suppose $h(x_1)=y_1$. Then $y_2\in F(h(x_1))=h(G(x_1))$ and this implies that there exists some $x_2\in G(x_1)$ such that $y_2=h(x_2).$ Similarly $y_3\in F(h(x_2))=h(G(x_2))$ and then there exists some $x_3\in G(x_2)$ with $y_3=h(x_3).$ In this process we can construct a forward orbit $\bar{x}=(x=x_1,x_2,\dots,x_n,\dots)$ with $y_i=h(x_i)$ for each $i.$ Let $U$ be an opene in $X.$ Then there exists an opene $V$ in $Y$ such that $h^{-1}(V)\subseteq U.$ Now there exists some $n\in \N,$ with $y_n\in V.$ Consequently, $x_n\in U.$ This proves that $\cO_G^{\oplus}(\bar{x})$ is dense in $X.$
	\end{enumerate}
\end{proof}

\begin{Theorem}
	Let $\{(X_i,F_i):i\in \mathbb{N}\}$ be a family of CR-dynamical systems together with irreducible epimorphisms $\{h_i:(X_{i+1},F_{i+1})\to (X_i,F_i):i\in \mathbb{N}\}$ and $(X,F)$ be the associated inverse limit. If there is some $m\in \N$ such that $(X_m,F_m)$ is point transitive (strongly point transitive) then $(X,F)$ has the corresponding property.
\end{Theorem}
\begin{proof}
	\begin{enumerate}
		\item Suppose that $(X_m,F_m)$ is point transitive. Since each $h_i$ is an irreducible epimorphism, from Theorems \ref{semiconju} and \ref{irrex}, we can conclude that  each $(X_i,F_i)$ is point transitive and so we can construct a point $\mathbf{x}=(x_k)\in X$ such that for each $k,$ $\displaystyle{\bigcup_{i=0}^{\infty}}F_k^{i}(x_k)$ is dense in $X_k.$ Let $U$ be opene in $X.$ Then from Lemma \ref{inverselimit}, there exists $n\in \N$ and opene $U_n\in X_n$ such that $p_n^{-1}(U_n)\subseteq U.$ Furthermore, there exists some $r\in \N$ with $F_n^r(x_n)\cap U_n = F_n^r(p_n(\mathbf{x}))\cap U_n\neq \emptyset,$ which implies that $p_n(F^r(\mathbf{x}))\cap U_n \neq \emptyset.$ Hence $F^r(\mathbf{x})\cap U \neq \emptyset.$ Therefore $(X,F)$ is point transitive.   
		\item Suppose that $(X_m,F_m)$ is strongly point transitive. Since each $h_i$ is an irreducible epimorphism, from Theorem \ref{semiconju} and \ref{irrex}, we can conclude that  each $(X_i,F_i)$ is strongly point transitive and also we can construct a point $\mathbf{x}=(x_k)\in X$ and a forward orbit $\bar{\mathbf{x}}=(\mathbf{x}=\mathbf{x}_1,\mathbf{x}_2,\dots,\mathbf{x}_i,\dots)\in T_F^+(\mathbf{x}),$ where for each $i\in \N,$ $\mathbf{x}_i=(x^i_k).$ It is evident that for each $k\in \N,$ $\bar{x_k}=(x_k=x^1_k,x^2_k,\dots)\in T_{F_k}^+(x_k)$ and $\cO_{F_k}^{\oplus}(\bar{x_k})$ is dense in $X_k.$ Let $U$ be opene in $X.$ Then from Lemma \ref{inverselimit}, there exists $n\in \N$ and opene $U_n\in X_n$ such that $p_n^{-1}(U_n)\subseteq U.$ Furthermore, there exists some $r\in \N$ with $x^r_n\in U_n \implies p_n(\mathbf{x}_r)\in U_n,$ which implies that $\mathbf{x}_r\in U.$ Therefore $(X,F)$ is strongly point transitive. 
	\end{enumerate}
\end{proof}

\bigskip

\begin{Definition}
	Let $(X,G)$ be a CR-dynamical system where $G$ is a suitable relation. Then  $(X,G)$ is called \emph{Suitably Topologically Transitive} or \emph{Suitably Transitive} if for all opene subsets $U,V \subseteq X$, there exists $n \in \N$
	such that $G^{\bullet n}(U) \cap V\neq \emptyset$.\par 
	 \end{Definition}
 
 \begin{Remark}
 	$(X,G,\bullet)$ is topologically transitive $ \iff $ $(X,G)$ is suitably transitive.
 \end{Remark}

Generally the concept of topological transitivity and point transitivity are mutually independent in a CR-dynamical system. But when we consider suitable dynamics we can relate suitable transitivity with point transitivity.

\bigskip

We now discuss the relation between transitivity of a closed relation with its selection functions.

\begin{Theorem}\label{4}
	Let $(X,F)$ be a CR-dynamical system with $F$ being a suitable relation on $X$ and $f$ be any selection function of $F$. Then the following conditions are equivalent:\par 
	\begin{enumerate}
		\item $(X,F)$ is suitably transitive.
		\item For any two opene subsets $U$ and $V$ of $X$ there exists $n\in \N$ such that $(U\cap f^{-n}(V))^{\circ} \neq \emptyset$.
		\item For any two opene subsets $U$ and $V$ of $X$ there exists $n\in \mathbb{N}$ such that $(U\cap F^{-n}(V))^{\circ} \neq \emptyset$.
		\item For every opene $U\subseteq X,$ $\displaystyle{\bigcup_{k=1}^{\infty}}(F^{-k}(U))^{\circ}$ is dense in $X.$
		\item For every opene subsets $U$ and $V$ of $X,$ the set $\{x\in U:$ $\mathcal{O}^{+}_F(x)\cap V \neq \emptyset\}$ is of the second category. 
		\item The set $\{x\in X: \overline{\mathcal{O}^{\oplus}_f(x)}=X\}$ is residual.  
		\item The set $\{x\in X: \overline{\mathcal{O}^{+}_F(x)}=X\}$ is residual. 
		\item There exists $x_0\in X$ such that $\mathcal{O}^{\oplus}_f(x_0)$ is dense in $X$ and $\mathcal{O}^{\oplus}_f(x_0)\subseteq C_f$.
		
	\end{enumerate}
	
\end{Theorem}
\begin{proof}
	$(1)\implies (2)$ is clear from Lemma \ref{suit int} and $(2)\implies (3)\implies (4)$ and $(6)\implies (7)$ are obvious.\par 

$(4)\implies (5)$ Let $U$ and $V$ be opene in $X.$ Then there exists $n\in N$ such that $(U\cap F^{-n}(V))^{\circ} \neq \emptyset.$ Clearly, $W=(U\cap F^{-n}(V))^{\circ}$ is an opene set and $W\subseteq \{x\in U: \mathcal{O}^{+}_G(x)\cap V \neq \emptyset\}.$ Being an opene set, $W$ is of the second category, which completes the proof.
	
	$(5)\implies (6)$ Let $\{V_1,V_2,V_3,\dots \}$ be a countable base for $X$. For each $n\in \N,$ let us consider the set $\cO_F^n=\{x\in X:$ $\mathcal{O}^{+}_G(x)\cap V_n \neq \emptyset\}$. It is obvious that each $\cO_F^n$ is of second category. Then from Corollary \ref{3}, we can conclude that $\mathcal{O}_F^n\cap C_f^{\infty}\neq \emptyset$, otherwise $\mathcal{O}_F^n$ will be of first category. Let $y\in \mathcal{O}_F^n \cap C_f^{\infty}$. Then clearly, $F^k(y)=f^k(y)$ as $f$ is continuous at $f^k(y)$ for every $k\in \mathbb{N}$. Furthermore we can obtain $n_0\in \mathbb{N}$ such that $f^{n_0}(y)\in V_n$. Since $f^{n_0}$ is continuous at $y,$ there exists an opene $W_y$ containing $y$ such that $f^{n_0}(W_y)\subseteq V_n$ and which implies that $W_y\subseteq \mathcal{O}_F^n$. Let us take $U_n=\displaystyle{\bigcup_{y\in \mathcal{O}_F^n\cap C_f^{\infty}}}W_y$. Clearly, $U_n$ is opene with $\mathcal{O}_F^n\cap C_f^{\infty}\subseteq U_n \subseteq \mathcal{O}_F^n$. Moreover, $U_n$ is dense in $X$. Indeed for each $m\in \mathbb{N}$, the set $\{x\in V_m:\mathcal{O}^{+}_F(x)\cap V_n \neq \emptyset\}\cap C_f^{\infty}\neq \emptyset$. Consequently, $\mathcal{O}_F^n\cap C_f^{\infty}\cap V_m \neq \emptyset$ and so $U_n\cap V_m \neq \emptyset$ for every $m\in \mathbb{N}$. Put $A_0=(\displaystyle{\bigcap_{n=1}^{\infty}}U_n)\cap C_f^{\infty}$. Then $A_0$ is residual and for every $x\in A_0$, $\mathcal{O}_f^{+}(x)$ is dense in $X$.\par 
	$(7)\implies (8)$ Since $C_f^{\infty}$ is residual, we have $\{x\in X: \overline{\mathcal{O}^{+}_F(x)}=X\}\cap C_f^{\infty}\neq \emptyset$. Hence there exists $x_0\in X$ such that $\mathcal{O}^{\oplus}_f(x_0)$ is dense in $X$ and $\mathcal{O}^{\oplus}_f(x_0)\subseteq C_f$.\par 
	$(8)\implies (1)$ Let $U$ and $V$ be opene subsets of $X$. Then from the given condition we can obtain $k\in \mathbb{N}$ such that $f^k(x_0)\in U$. Take $W=V\setminus \{f^i(x_0):i=0,1,\dots k\}$. Clearly, $W$ is also an opene set. Consequently, there is $s>k$ such that $f^s(x_0)\in W$. Since $f^n(x_0)\in C_f$ for every $n\in \N$, $(f^k(x_0),f^s(x_0))\in F^{\bullet (s-k)}$. Hence $(X,F)$ is suitably transitive.
\end{proof}
\begin{Corollary}\label{qc trans}
	Let $(X,f)$ be a quasi-continuous dynamical system and $f$ be a quopen map on $X$. Then the following conditions are equivalent:\par 
	$(1)$ $(X,\overline{f})$ is suitably transitive.\par 
	$(2)$ The set $\{x\in X: \overline{\mathcal{O}^{\oplus}_f(x)}=X\}$ is residual.\par
	$(3)$ $(X,f)$ is transitive. 
\end{Corollary}
\begin{proof}
	The proof easily follows from Theorem \ref{quasi dyn}.
\end{proof}

\begin{Remark}
 The condition $(5)$ of Theorem \ref{4} was called `strongly transitive' in \cite{multifunction}. They had considered cm functions, whose closure need not be  suitable relations. Thus our result can be considered to be more general in context of  relations, however it is true in some other cases too. \end{Remark}

Every quasi-continuous system $(X,f),$ where $f$ is a quopen map induces a CR-dynamical sytem $(X,F),$ with $F$ a suitable relation. Moreover, from Corollary \ref{qc trans} we can conclude that transitivity of $(X,f)$ is equivalent to the suitable transitivity of $(X,F).$ Now we present some examples which illustrate difference between transitivity of CR-dynamical system, suitable transitivity and transitivity of a quasi-continuous system. We note that suitable transitivity is a generalization of transitivity of both quasi-continuous maps as well as the continuous maps. 

\begin{Example} \label{nointerior}
	 Let $X=[0,1]$ and a closed relation $F$ be defined as $$F=\displaystyle{\bigcup_{n=1}^{\infty}}\big(\big[\frac{1}{n+1},\frac{1}{n}\big]\times \big[\frac{1}{n+1},\frac{1}{n}\big]\big)\bigcup \big\{(0,0)\big\}.$$
	  Let $U$ and $V$ be opene in $X.$ Then there exist $m,n\in N$ such that $[\frac{1}{(m+1)},\frac{1}{m}]\cap U\neq \emptyset$ and $[\frac{1}{(n+1)},\frac{1}{n}]\cap V\neq \emptyset.$ Then evidently $F^{(|m-n|+1)}(U)\cap V \neq \emptyset.$ Hence $(X,F)$ is topologically transitive, but $F$ is not a suitable relation. So we can not define suitable transitivity for $(X,F).$
		
\end{Example}

\begin{Example} \label{longtent}
	Suppose that $X=[0,1]$ and $f:X\to X$ is a function with $f(0)=0,$ $f(\frac{1}{2n-1})=0$ and $f(\frac{1}{2n})=1$ for every $n\in \N$ and $f$ is linear elsewhere. Take $F=\overline{f}.$ Clearly, $F$ is suitable and for every selection function $g$ we have $F^{\bullet n}=\overline{g^n}=F^n$ for each $n\in\N.$ Then from Theorem \ref{qc}, we can conclude that $(X,g)$ is a quasi-continuous system. Moreover, for every opene $U,$ there exists $n\in \N$ such that $F^n(U)=[0,1]$. Hence $(X,F)$ is suitably transitive and also from Corollary \ref{qc trans}, each $(X,g)$ is transitive.
\end{Example}
\begin{Example}\label{tent}                                 Let $T$ be the tent map defined as $T(x)=2x$ for $x\in [0,\frac{1}{2}]$ and $T(x)=2-2x$ for $x\in [\frac{1}{2},1].$ Suppose that $X=[-1,1]$ and $f:X\to X$ is a function such that $f(x)=T(x)$ for every $x\in [0,1]$ and $f(-\frac{1}{2n-1})=-1,$ $f(-\frac{1}{2n})=1$ for every $n\in \N$ and $f$ is linear elsewhere. Take $F=\overline{f}.$ Clearly, $F$ is suitable and for every selection function $g$ we have $F^{\bullet n}=\overline{g^n}$ for each $n\in\N.$ Then from Theorem \ref{qc}, we can conclude that $(X,g)$ is a quasi-continuous system. But $F\circ F \neq F\bullet F.$ Here $F\circ F=(F\bullet F) \cup (\{1\}\times [-1,1]).$ \par 

For every opene $U \subseteq X$ there exists $n\in \N,$ such that $0\in F^n(U).$ Consequently, $F^{(n+1)}(U)=[0,1].$ Hence $(X,F)$ is transitive. On the other hand if we take $U=(\frac{1}{2},1)$ and $V=(-1,-\frac{1}{2}).$ Then $\{x\in U:$ $\mathcal{O}^{+}_F(x)\cap V \neq \emptyset\}\subseteq \{\frac{m}{2^n}:m,n\in \N\},$ which is not of second category. Hence from Theorem \ref{4}, we can conclude that $(X,F)$ is not suitably transitive. Also the quasi-continuous system $(X,g)$ is not transitive for every selection function $g$. This example show that for suitable $G$, the CR-dynamical system $(X,G)$ and  its suitable version $(X,G, \bullet)$ are not conjugate and give different dynamics. 
\end{Example}

\begin{Example}\label{tent}                                
	Let $T$ be the tent map defined as $T(x)=2x$ for $x\in [0,\frac{1}{2}]$ and $T(x)=2-2x$ for $x\in [\frac{1}{2},1].$ Suppose that $X=[-1,1]$ and $f:X\to X$ is a function such that $f(x)=T(x)$ for every $x\in [0,1]$ and $f(-\frac{1}{2n-1})=-1,$ $f(-\frac{1}{2n})=1$ for every $n\in \N$ and $f$ is linear elsewhere. Take $F=\overline{f}.$ Clearly, $F$ is suitable and $F^{\bullet n}=\overline{f^n}$ for each $n\in\N.$ Then from Theorem \ref{qc}, we can conclude that $(X,f)$ is a quasi-continuous system. But $F\circ F \neq F\bullet F.$ Here $F\circ F=(F\bullet F) \cup (\{1\}\times [-1,1]).$ Note that if we consider a selection function $g$ such that $g(x)=f(x)$ for $x\neq 0$ and $g(0)=1,$ then $g^2$ is not quasi-continuous. \par
	
	For every opene $U \subseteq X$ there exists $n\in \N,$ such that $0\in F^n(U).$ Consequently, $F^{(n+1)}(U)=[-1,1].$ Hence $(X,F)$ is transitive. On the other hand if we take $U=(\frac{1}{2},1)$ and $V=(-1,-\frac{1}{2}).$ Then $\{x\in U:$ $\mathcal{O}^{+}_F(x)\cap V \neq \emptyset\}\subseteq \{\frac{m}{2^n}:m,n\in \N\},$ which is not of second category. Hence from Theorem \ref{4}, we can conclude that $(X,F)$ is not suitably transitive. Also, the quasi-continuous system $(X,f)$ is not transitive. This example shows that for suitable $G$, the CR-dynamical system $(X,G)$ and its suitable version $(X,G, \bullet)$ are not conjugate and give different dynamics.
\end{Example}

\bigskip 

Thus we have, 

\bigskip

	\begin{tabular}{|c|}
	\hline\\
	
	point transitive  $\Longleftarrow$   strongly point transitive $\Longrightarrow$ topologically transitive.\\

	suitably transitive $\Longrightarrow$  topologically transitive.\\
	
	\hline
\end{tabular}

\bigskip

\subsection{Weakly Mixing,  Mixing and  suitable forms} 	Let $(X,G)$ be a CR-dynamical system. We recall that $G$ induces the closed relation
$$G  \times G \subseteq X  \times X \times X  \times X \ \text{by} \ \{(x_1, y_1, x_2, y_2) : (x_1, x_2),(y_1, y_2) \in G\}.$$

$ G \times G $ is just the set product of the two relations with the second and third coordinates switched.

We can inductively define, $G^{[n]} = \underbrace{G  \times \ldots \times G}_{n \ \text{times}} \subseteq X^{2n}  $  as 
$$G^{[n]} = \{(x_1, y_1, x_2, y_2, \ldots, x_n,  y_n) : (x_1, y_1),(x_2, y_2) , \ldots, (x_n,y_n)\in G\}.$$

\begin{Definition}
	$(X,G)$ is called 
	
	- \emph{Mixing or Strongly Mixing } if for every pair of opene sets $U,V\subseteq X$, there exists $N\in \mathbb{N}$, such that $G^n(U)\cap V \neq \emptyset$ for all $n\geq N$.
	
	- \emph{Weakly Mixing} if the product system $(X\times X, G\times G)$ is topologically transitive.\par 
	
\end{Definition}

Note that  mixing $\Longrightarrow$   weakly mixing $\Longrightarrow$  topologically transitive.

\bigskip

Note that Example \ref{nointerior} is both weakly mixing and  mixing. We consider its variation here.

\begin{Example} \label{everything}
	Let $X=[0,1]$ and a closed relation 
	$$G=\big(\displaystyle{\bigcup_{n=1}^{\infty}}[\frac{1}{2^{n}},\frac{1}{2^{n-1}}]\times [\frac{1}{2^{n}},\frac{1}{2^{n-1}}] \big) \cup \big(\displaystyle{\bigcup_{n=0}^{\infty}}[\frac{1}{3.2^{n}},\frac{1}{3.2^{n-1}}]\times [\frac{1}{3.2^{n}},\frac{1}{3.2^{n-1}}] \big)\cup\{(0,0)\}.$$ Then $(X,G)$ is both weakly mixing and mixing.

\end{Example}

We skip the trivial proof of the following:

\begin{Theorem}  For a CR-system $(X,G)$, the
	following are equivalent.
	\begin{enumerate}
		\item  The system is  mixing.

		\item  For every pair of opene
		sets $U, V \subseteq X$ the set $N(U,V)$ is cofinite.

	\end{enumerate}
		
\end{Theorem}

\begin{Definition}
	Let $(X,G)$ be a CR-dynamical system where $G$ is a suitable relation. Then  
	$(X,G)$ is called 
	
	- \emph{Suitably Mixing}   if for every pair of opene sets $U,V\subseteq X$, there exists $N\in \mathbb{N}$, such that $G^{\bullet n}(U)\cap V \neq \emptyset$ for all $n\geq N$.
	
	- \emph{Suitably Weakly Mixing} if the product system $(X\times X, G\times G)$ is suitably topologically transitive, i.e., if for pairs of opene subsets $U_1, U_2 ,V_1, V_2 \subseteq X$, there exists $n \in \N$
	such that $G^{\bullet n}(U_i) \cap V_i \neq \emptyset$ for $i=1,2$.
	
\end{Definition}

\begin{Remark}
	$(X,G)$ is suitably weakly mixing or suitably mixing is equivalent to saying that $(X,G, \bullet)$ is weakly mixing or mixing.
	
	Clearly, in dynamics of continuous maps the concepts of weakly mixing and suitably weakly mixing or mixing and suitably mixing are same.
\end{Remark}

Example \ref{longtent} is both suitably weakly mixing and suitably mixing. However, Examples \ref{nointerior} and \ref{everything} are weakly mixing and mixing but being non suitable relations cannot be of a suitable form. Moreover, Example \ref{tent} is weakly mixing as well as mixing, where the given relation is suitable but it is not suitably weakly mixing and so not suitably mixing.

\bigskip

Note that  suitably mixing $\Longrightarrow$   suitably weakly mixing $\Longrightarrow$  suitably transitive. 

\bigskip

We define the following notation for the CR-dynamical system $(X,G, \bullet)$;
$$N^{\bullet}(U,V)= \{n\in \N: G^{\bullet n}(U)\cap V\neq \emptyset\}$$

For opene $U,V \subset X,$  from Lemma \ref{suit int}, we can say that $N^{\bullet}(U,V)=\{n\in \N: (U\cap G^{-n}(V))^{\circ}\neq \emptyset\}$

Similarly for $x\in X,$ and opene $U \subset X$ one can define 
$$N^\bullet(x,U)=\{n\in \N: G^{\bullet n}(x)\cap U\neq \emptyset\}$$  $$N^\bullet(U,x)=\{n\in \N: x\in G^{\bullet n}(U)\}$$

\begin{Theorem} \label{wm} Let $G$ be a suitable closed relation on $X$ and $g$ be a selection function of $G$. Then for the CR-dynamical system  $(X,G)$  the
	followings are equivalent.
	\begin{enumerate}
		\item  $(X,G)$ is suitably weakly mixing.
		
		\item For any two  opene sets $U, V$ in $X$, there exists $N \in \N$ such  that $G^{\bullet N}(U) \cap V \neq \emptyset$ and $G^{\bullet N}(U) \cap U \neq \emptyset$.

		\item For any two  opene sets $U, V$ in $X$, there
		exists $N \in \N$ such  that $(U \cap g^{-N}(V))^{\circ} \neq
		\emptyset$ and $(U\cap g^{-N}(U))^{\circ} \neq \emptyset$.
		
		\item $N^{\bullet}(U,V) \cap N^{\bullet}(U,U)$ $ \neq  \emptyset$ for every pair of opene  $U,V \subseteq X$.
		
		\item For every $N \in \N,$ the product system $(X^N, G^{[N]})$ is suitably transitive.

	\end{enumerate}
	
\end{Theorem}
\begin{proof} We note that $(2) \Longleftrightarrow (3) \Longleftrightarrow (4)  $ is follows from Lemma \ref{suit int}.
	
	We see that	$(1) \Longrightarrow (2)$ holds from the definition, and $(5) \implies (1)$ is vacuous. 
	
	We prove $(3) \Longrightarrow (1)$. This is just  characterization of suitable weak mixing relations similar to the one given by Petersen for continuous maps in \cite{pk}. We prove the same for suitable CR-dyanmical systems. Consider opene sets $U_1,U_2,V_1,V_2 \subseteq X$. From given condition we have

	Let	$ U = (g^{-k_1}(U_1) \cap V_1)^{\circ} \neq \emptyset $ and 		
		$ V = (g^{-k_2}(U)\cap (g^{-k_1}(U_2))^{\circ})^{\circ} \neq \emptyset $,
	
	for some $k_1, k_2 \in \N$. Take $k_3 \in \N$ such that
		
		$$ (g^{-k_3}(V) \cap V)^{\circ} \neq \emptyset \ \text{and} \ (g^{-k_3}(V) \cap V_2)^{\circ} \neq \emptyset .$$

Now $g^{-k_1}(g^{-(k_2+k_3)}(U_1) \cap U_2) = g^{-(k_1+k_2+k_3)}(U_1) \cap g^{-k_1}(U_2)$
	
$\supseteq g^{-(k_1+k_2+k_3)}(U_1) \cap g^{-(k_2+k_3)}(V_1) \cap g^{-k_1}(U_2)$ 
	
$ = g^{-(k_2+k_3)}(g^{-k_1}(U_1) \cap V_1) \cap g^{-k_1}(U_2)$
	
$ \supseteq g^{-(k_2+k_3)}(U) \cap g^{-k_1}(U_2) \supseteq g^{-k_3}(V) \cap g^{-k_1}(U_2)$
	
	$ \supseteq g^{-k_3}(V) \cap V \neq \emptyset$. Thus,

	$$	W=(g^{-k_1}(g^{-(k_2+k_3)}(U_1) \cap U_2))^{\circ} \neq \emptyset.$$
	
Then from Corollary \ref{3}, there exists a point $x\in W\cap C_g^{\infty},$ which implies that $g^{k_1+k_2+k_3}(x)\in U_1$ and $g^{k_1}(x)\in U_2.$ Since $g^{k_2+k_3}$ is continuous at $g^{k_1}(x),$ there is an opene set $W_0$ containing $g^{k_1}(x)$ such that $W_0\subseteq g^{-(k_2+k_3)}(U_1).$ Hence $(g^{-(k_2+k_3)}(U_1) \cap U_2)^{\circ} \neq \emptyset.$

	Also $g^{-(k_1+k_2+k_3)}(U_1) \cap g^{-(k_2+k_3)}(V_1) \cap V_2 = g^{-(k_2+k_3)}(g^{-k_1}(U_1) \cap V_1) \cap V_2 \supseteq g^{-(k_2+k_3)}(U) \cap V_2 \supseteq g^{-k_3}(V) \cap V_2 \neq \emptyset$. And so we can conclude that

		$$(g^{-(k_2+k_3)}(V_1) \cap V_2)^{\circ} \neq \emptyset.$$
		
		Then we get $(G^{-(k_2+k_3)}(U_1) \cap U_2)^{\circ} \neq \emptyset$ and $(G^{-(k_2+k_3)}(V_1) \cap V_2)^{\circ} \neq \emptyset$ i.e. $  U_1 \cap G^{\bullet (k_2+k_3)}(U_2) \neq \emptyset$ and $ V_1 \cap G^{\bullet (k_2+k_3)}(V_2) \neq \emptyset$.
	
	Hence by Theorem \ref{4}, $(X,G)$ is suitably weakly mixing.

	\bigskip
	
	The proof is complete if we show $(4) \Longrightarrow (5)$. But this is just a consequence of the \emph{Furstenberg Intersection Lemma} for CR-dynamical systems, and we see that the proof is similar to that for continuous maps.

	\emph{[Furstenberg Intersection Lemma] For a CR-dynamical system $(X,G, \bullet)$ if $N^{\bullet}(U,V) \cap N^{\bullet}(U,U)$ $ \neq  \emptyset$ for every pair of opene  $U,V \subseteq X$, then for all opene  $U_1,V_1,U_2,V_2 \subseteq X$ there exist opene  $U_3,V_3 \subseteq X$ such that
	$$ N^{\bullet}(U_3,V_3) \quad \subseteq \quad  N^{\bullet}(U_1,V_1) \cap  N^{\bullet}(U_2,V_2). $$}	
	
	We proceed to prove this statement. First of all since $(X,G)$ is transitive, $G$ is surjective. We see that since $ N^{\bullet}(U_1,V_1) \neq \emptyset$  there exists $n_1 \in\N$ such that $U_0 = (U_1 \cap G^{-n_1}(V_1))^{\circ}$ is nonempty.
	$ N^{\bullet}(U_0,U_2) \neq \emptyset$ implies there exists $n_2 \in\N$ such that $U = (U_1 \cap G^{-n_1}(V_1) \cap G^{-n_2}(U_2))^{\circ} \neq \emptyset$.
	
	From Lemma \ref{ni}, $G^{-n_1-n_2}(V_2)$ has non-empty interior. Then,

$$N^{\bullet}(U,U)\cap N^{\bullet}(U,G^{-n_1-n_2}(V_2)) \subseteq N^{\bullet}(U_1,G^{-n_2}(U_2))\cap N^{\bullet}(G^{-n_1}(V_1),G^{-n_1}(G^{-n_2}(V_2))) $$

One can observe that $k\in N^{\bullet}(G^{-n_1}(V_1),G^{-n_1}(G^{-n_2}(V_2)))  \implies G^{-n_1}(V_1)\cap G^{-(n_1+k+n_2)}(V_2)\cap C_g^{\infty}\neq \emptyset,$ which implies that there exists $x\in C_g^{\infty}$ such that $g^{n_1}(x)\in V_1$ and $g^{n_1+n_2+k}(x)\in V_2$. Moreover, $g^{n_2+k}$ is continuous at $g^{n_1}(x)$ and so there exists an opene set $W$ containing $g^{n_1}(x)$ such that $W\subseteq G^{-(k+n_2)}(V_2).$ Therefore $N^{\bullet}(G^{-n_1}(V_1),G^{-n_1}(G^{-n_2}(V_2)))\subseteq N^{\bullet}(V_1,G^{-n_2}(V_2)).$ Thus

		$$\emptyset \neq N^{\bullet}(U,U)\cap N^{\bullet}(U,G^{-n_1-n_2}(V_2))
		\subseteq   N^{\bullet}(U_1,G^{-n_2}(U_2))\cap N^{\bullet}(V_1,G^{-n_2}(V_2)).$$

	Fix $n_0 \in N^{\bullet}(U_1,G^{-n_2}(U_2))\cap N^{\bullet}(V_1,G^{-n_2}(V_2))$. With $n = n_0 + n_2$ consider the sets
	$U_3 = (U_1 \cap G^{-n}(U_2))^{\circ}, V_3 = (V_1 \cap G^{-n}(V_2))^{\circ}$.

	Let $k \in N^{\bullet}(U_3, V_3)$. Then $(G^{-k}(V_3) \cap U_3)^{\circ} \neq \emptyset$. 
	
	That
	is $(G^{-k}(V_1) \cap G^{-n-k}(V_2) \cap U_1 \cap G^{-n}(U_2))^{\circ} \neq \emptyset$.
	
	Hence $k \in N^{\bullet}(U_1, V_1) \cap N^{\bullet}(G^{-n}(U_2), G^{-n}(V_2)) \subseteq N^{\bullet}(U_1, V_1) \cap N^{\bullet}(U_2, V_2)$.

	Note that $N^{\bullet}(G^{-n}(U_2),G^{-n}(V_2)) \subseteq N^{\bullet}(U_2,V_2)$
	and so our assertion  follows.
	
	We now prove that $(5)$ holds. We prove this by induction on $\N$. Consider opene $U_1, U_2, \ldots U_k, U_{k+1}, V_1, V_2, \ldots, V_k, V_{k+1} \subseteq X$, and assume that $ G^{-n}(U_n) \cap V_n \neq \emptyset$ for $n = 1,2, \ldots, k$. 
	
	Then there exists opene sets $U, V \subseteq X$ such that $N^{\bullet}(U,V) \subseteq \bigcap \limits_{n=1}^k N^{\bullet}(U_n, V_n)$.
	
	But we have for opene $U,V,U_{k+1},V_{k+1}$, $N^{\bullet}(U,V) \cap N^{\bullet}(U_{k+1}, V_{k+1}) \neq \emptyset$. This proves our assertion.

\end{proof}

\begin{Remark}
	We note that it is possible for a weakly mixing CR-dynamical system $(X,G)$ that for every $N \in \N,$ the product system $(X^N, G^{[N]})$ is  transitive, even when $G$ is not suitable.
	
	Consider Example \ref{everything}. Here for every $N \in \N,$ the product system $(X^N, G^{[N]})$ is  transitive. 
\end{Remark}

\bigskip

\noindent\begin{tabular}{|c|}
	
	\hline\\
	
	\textbf{Question:} For a weakly mixing CR-dynamical system $(X,G)$, where $G$ need \\ not be	suitable,  is it always true that $(X^{[n]},G^{[n]})$ is transitive for every $n\in \N.$\\
	
	\hline
\end{tabular}

\bigskip 

Recall that that a subset $M \subseteq \mathbb{N}$ is called \emph{thick} if it contains arbitrary long segments of consecutive elements of $\mathbb{N}$, i.e. for given $n\in \mathbb{N}$, there exists $k\in \mathbb{N}$ such that $k,k+1,\dots,k+n\in M$.

\begin{Theorem} Suppose that $G$ is a suitable closed relation on $X$ and $(X,G)$ is a suitably weakly mixing CR-dynamical system. Then for every pair of opene sets $U,V\subseteq X$, $N^{\bullet}(U,V)$ is a thick subset of $\mathbb{N}$.
\end{Theorem}
\begin{proof} Let $(X,G)$ be a suitably weakly mixing CR-dynamical system and $U,V$ be two opene subsets of $X$. Then from Lemma \ref{ni}, it is clear that for every $n\in \mathbb{N}$, $(G^{-n}(V))^{\circ}\neq \emptyset$. Now from Theorem \ref{wm}, we can conclude that for every $n\in \mathbb{N}$ there exists $k\in N^{\bullet}(U,V)\cap N^{\bullet}(U,G^{-1}(V))\cap \dots \cap N^{\bullet}(U,G^{-n}(V))$, which implies that $(U\cap G^{-k}(V))^{\circ}\neq \emptyset$, $(U\cap G^{-(k+1)}(V))^{\circ}\neq \emptyset,\dots, (U\cap G^{-(k+n)}(V))^{\circ}\neq \emptyset$. Hence $k,(k+1),\dots,(k+n) \in N^{\bullet}(U,V)$ and so $N^{\bullet}(U,V)$ is a thick subset of $\mathbb{N}$. 
	
\end{proof}

We now focus on weakly mixing property of a quasi-continuous dynamical system.
\begin{Definition}
	A quasi-continuous dynamical system $(X,f)$ is called weakly mixing if $(X\times X,f\times f)$ is topologically transitive.
\end{Definition}
Clearly, $(X,f)$ is weakly mixing if and only if  for opene sets $U_1,U_2,V_1,V_2$ we have $N(U_1,V_1)\cap N(U_2,V_2)\neq \emptyset$. 

\bigskip

The following result makes a connection between weakly mixing quasi-continuous system and suitably weakly mixing.
\begin{Theorem}
	Let $(X,f)$ be a quasi-continuous system with $f$ quopen. Then the followings are equivalent:
\begin{enumerate}
	\item $(X,\overline{f})$ is suitably weakly mixing.
	\item $(X,f)$ is weakly mixing.
	\item For any two  opene sets $U, V$ in $X$, there exists $N \in \N$ such  that $(U \cap f^{-N}(V))^{\circ} \neq \emptyset$ and $(U\cap f^{-N}(U))^{\circ} \neq \emptyset$.
	\item For every $N \in \N,$ the product system $(X^N, f^{(N)})$ is transitive.	
\end{enumerate}
\end{Theorem}
\begin{proof}
	The proof  follows from Theorem \ref{quasi dyn}.
\end{proof}

\bigskip

Thus, for  $(X,G)$ we have:

\bigskip
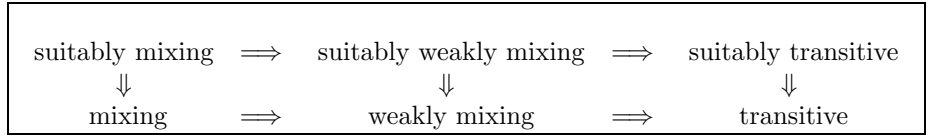
\begin{tabular}{|c|}
	\hline\\
	
	{{\small
			$\begin{array}{ccccc}
				\text{suitably mixing}  & \Longrightarrow & \text{ suitably weakly mixing} & \Longrightarrow &  \text{ suitably transitive}\\
				\Downarrow &    & \Downarrow &   & \Downarrow\\
				\text{ mixing}  & \Longrightarrow & \text{  weakly mixing} & \Longrightarrow &  \text{  transitive}
			\end{array}$}}	\\
	
	\hline
\end{tabular}

\bigskip

\bigskip

	\subsection{Minimal and suitably minimal}
	
	Some of our definitions and results overlap with \cite{mini} and some results of \cite{mini} are presented here in much simpler form.

\begin{Definition} Let $(X,G)$ be CR-dynamical system. Then
	\begin{enumerate}
		\item $(X,G)$ is called \emph{minimal} if $\overline{\cO^+_G(x)}=X$, for each $x\in X$. (It is called $3^{\oplus}$-minimal in \cite{mini}).
	\item $(X,G)$ is called \emph{strongly minimal} if for every $x\in X$ and evey forward orbit $\bar{x}\in T_G^+(x),$ $\mathcal{O}_G^{\oplus}(\bar{x})$ is dense in $X$.(It is called $1^{\oplus}$ minimal in \cite{mini}).
	\end{enumerate}
\end{Definition}

Note that strongly minimal $\Longrightarrow$ minimal.	And the converse need not hold even for suitable relations.

\begin{Example} \label{irr}
		Consider  $X = \mathbb{S}$,	the unit circle bijective with $[0,1)$ and let $ \alpha \in (0,1) $ be  irrational.
		
	Let $ T(x)= x+\alpha $ be the irrational rotation on $X$.  
	Define a closed relation $ G  \subseteq X \times X$  as: 
	
$$ 	G=\big\{\big(x, T(x)\big): x\in \big[0,\frac{1}{2}\big]\cup \big[\frac{1}{2}+\frac{\alpha}{2}, 1 -\frac{\alpha}{2}\big]\big\} \bigcup \big\{\big(x, x+\frac{1}{2}\big) : x\in \big[\frac{1}{2}, \frac{1}{2}+\frac{\alpha}{2}\big] \cup \big[1 - \frac{\alpha}{2}, 1\big]\big \} $$
	
	That means $ G $ is irrational rotation everywhere except $(\frac{1}{2}, \frac{1}{2}+\frac{\alpha}{2})\cup (1 - \frac{\alpha}{2}, 1)$ where it is a rotation by $\frac{1}{2}$, and is multivalued on $\frac{1}{2}, \frac{1}{2} + \frac{\alpha}{2}, 1 - \frac{\alpha}{2}, 1$. 
	
	We note that $ G $ is a suitable relation, since any selection function of $ G $ is quasi continuous. 
	
	Now $G$ is minimal since $\overline{\cO^+_G(x)}=X$, for each $x \in X$.
	
	 However, for $x = \frac{1}{2}$, $T_G^+(\frac{1}{2}) \ni \bar{a}$ with $\bar{a} = \frac{1}{2}, 1, \frac{1}{2}, 1, \frac{1}{2}, 1, \ldots$ and $\mathcal{O}_G^{\oplus}(\bar{a}) = \{\frac{1}{2}, 1\}$. 	
	Hence $ G $ is not strongly minimal.

\end{Example}

\begin{Definition}  Let $(X,G)$ be a CR-dynamical system and $A\subseteq X.$
	\begin{enumerate}
		\item $A$ is called a \emph{minimal set} if $A$ is a closed weakly invariant subset of $(X,G)$ and the system $(A,G_A)$ is minimal.
		\item $A$ is called \emph{strongly minimal} if $A$ is a closed weakly invariant subset of $(X,G)$ and the system $(A,G_A)$ is strongly minimal.
	\end{enumerate}
\end{Definition}

\begin{Theorem}\label{minimality}
	Suppose that $(X,G)$ is a CR-dynamical system. Then the following conditions are equivalent.\par 
\begin{enumerate}
	\item $(X,G)$ is minimal.\par
	\item For every opene set $U\subseteq X$ and for every point $x\in X$, there exists $n\in \mathbb{N}$ such that $G^{n}(x)\cap U\neq \emptyset$.\par 
	\item For every $x\in X$, $\displaystyle{\bigcup_{n=1}^{\infty}}G^{n}(x)$ is dense in $X$.\par 
	\item For every opene set $U\subseteq X$ and every point $x\in X$, the set $N(x,U)$ is nonempty.\par 
	\item For every opene set $U\subseteq X$ and every point $x\in X$, the set $N(x,U)$ is infinite.\par
	\item If $E$ is nonempty and $+$invariant, then $E$ is dense in $X$. \par
	\item For every $x\in X,$ there exists $\bar{x}\in T_G^+(x)$ such that $\cO_G^{\oplus}(\bar{x})$ is dense in $X$ (This condition is called $2^{\oplus}$-minimality in \cite{mini}.)\par
	\item For every $x\in X$ there exists $\bar{x}\in T_G^+(x)$ such that for every opene $U\subseteq X,$ $N_{\bar{x}}(x,U)$ is infinite.
	\end{enumerate}

In particular, $(X,G)$ is  strongly  point transitive.
	
\end{Theorem}

\begin{proof}
	Since $X$ has no isolated point, the equivalence of $(1),(2), (3), (4)$ is trivial. $(8)\Longrightarrow (5)\Longrightarrow (4)$ is clear and $(3)\iff(6)$ is evident from the fact that for every $x\in E$, we have $\displaystyle{\bigcup_{n=1}^{\infty}}G^{n}(x)\subseteq E$.
	
	\medskip
	
	 We only prove $(2)\Longrightarrow (7) \Longrightarrow (8).$\par   
 
 For $(2)\Longrightarrow (7)$ suppose that $(X,G)$ is minimal.  Let $\{U_n\}$ be a countable base for $X$ and $x$ be any point in $X$. Then there there is some $n_1\in \N$ such that $G^{n_1}(x)\cap U_1\neq \emptyset.$ Let $y_1\in G^{n_1}(x)\cap U_1.$ Again there is some $n_2\in \N$ with $G^{n_2}(y_1)\cap U_2\neq \emptyset.$ Let $y_2\in G^{n_2}(y_1)\cap U_2.$ Proceeding in this way we have $y_k\in G^{n_k}(y_{k-1})\cap U_k$ for each $k>1.$ Then we can construct a forward orbit $\bar{x}=(x,\dots,y_1,\dots,y_2,\dots,y_k,\dots)\in T_G^+(x),$ where $y_1\in G^{n_1}(x)$ and $y_k\in G^{n_1+\dots+n_k}(x)$ for every $k>1.$ Clearly, $\cO_G^{\oplus}(\bar{x})$ is dense in $X.$\par
 
 \bigskip
 
	For $(7)\Longrightarrow (8)$ let $x\in X$ and $U$ be opene in $X.$ Suppose that $\bar{x}=(x,x_1,x_2,\dots)\in T_G^+(x)$ such that $\cO_G^{\oplus}(\bar{x})$ is dense in $X$. Assume that $N_{\bar{x}}(x,U)$ is finite. Let $N=\max\{n:n\in N_{\bar{x}}(x,U)\}.$ Now take $W=U\setminus \{x,x_1,\dots,x_N\}.$ Then there exists $M>N$ such that $x_M\in W,$ which is a contradiction. Hence $N_{\bar{x}}(x,U)$ is infinite.
	
\end{proof}

\textbf{Observation:} From $(2)\Longleftrightarrow (7)$ of the above theorem, we can conclude that the concepts of $2^{\oplus}$-minimality and $3^{\oplus}$-minimality of \cite{mini} are equivalent. Hence this  answers the problem 4.13 in \cite{mini}. 

\bigskip

We recall that a subset $M \subseteq \mathbb{N}$ is called \emph{syndetic} if there exists $N\in \mathbb{N}$ such that every interval of length $N$ in $\mathbb{N}$ meets $M$, i.e. given $k \in \N$ we have $\{ k+1, \ldots, k+N \} \cap M \neq \emptyset$.

\begin{Theorem}
	Let $(X,G)$ be a CR-dynamical system. Then the following are equivalent:
	\begin{enumerate}
		\item $(X,G)$ is strongly minimal.
		\item For every opene $U\subseteq X$ and every point $x\in X$ and every $\bar{x}\in T_G^+(x),$ the set $N_{\bar{x}}(x,U)$ is infinite.
		\item For every opene $U\subseteq X$ and every point $x\in X$ and every $\bar{x}\in T_G^+(x),$ the set $N_{\bar{x}}(x,U)$ is syndetic.
		\item If $E$ is nonempty, closed and weakly invariant then $E=X$. 
		\item If $E$ is nonempty and weakly invariant then $E$ is dense in $X.$
		\item For every point $x\in X$ and every $\bar{x}\in T_G^+(x)$ the set $\omega(\bar{x})=X.$ 
	\end{enumerate}
\end{Theorem}

\begin{proof}
	We see that $(1)\Longrightarrow (2)$ is same as $(7)\Longrightarrow (8)$ of Theorem \ref{minimality}.\par 
	For $(2)\Longrightarrow (3)$ we suppose that there exist opene $U\subseteq X$ and a forward orbit $\bar{x}\in T_G^+(x)$ such that $N_{\bar{x}}(x,U)$ is not syndetic. Since $N_{\bar{x}}(x,U)$ is infinite,  for all $k\in \N$ we can choose positive integers $n_1<n_2<\dots$ such that $\pi_{n_k}(\bar{x})\in U$ but $\{\pi_{n_k+m}(\bar{x}):m=1,\dots,k\}\cap U=\emptyset$.  Now from compactness, $(\pi_{n_k}(\bar{x}))$ has a cluster point say $y$. Without any loss of generality we assume that $(\pi_{n_k}(\bar{x})) \to y$. Then if necessary by subsequently passing to a subsequence, we have $(\pi_{n_k+1}(\bar{x})) \to y_1$, $(\pi_{n_k+2}(\bar{x})) \to y_2$, \ldots, $(\pi_{n_k+r}(\bar{x})) \to y_r$, for every $r\in \N$. Since $G$ is closed we have $(y,y_1)\in G$ and subsequently $(y_r,y_{r+1})\in G$ for every $r\in \N$. Thus we get the forward orbit $\bar{y}=(y,y_1,y_2,\dots)\in T_G^+(y)$. Since $(X,G)$ is strongly minimal there is some $l\in \N$ such that $y_l\in U$. Hence for some $N>l$ we have $\pi_{n_N+l}(\bar{x})\in U$, contradicting our assumption. This completes the proof.\par 
	
	\bigskip
	
	For $(3)\Longrightarrow (4)$ let $x\in E$. Since $E$ is closed and weakly invariant, there is a forward orbit $\bar{x}\in T_G^+(x)$ such that $\overline{\cO_G^{\oplus}(\bar{x})}\subseteq E$. For any opene $U \subseteq X$, $N_{\bar{x}}(x,U)$ is syndetic, and so $\cO_G^{\oplus}(\bar{x})$ is dense in $X$. Hence $E=X$.\par 
	
	\bigskip
	
	For $(4)\Longrightarrow (5)$ we observe from  Lemma \ref{cl of wi},  that $\overline{E}$ is also nonempty closed weakly invariant subset of $X$. Thus $\overline{E}$ and hence $E$ is dense in $X$.\par
	
	\bigskip
	
	For $(5)\Longrightarrow (6)$ let $x\in X$ and $\bar{x}\in T_G^+(x)$. From Lemma \ref{cl of wi}, it is clear that $\omega(\bar{x})$ is a nonempty closed weakly invariant subset of $X$. Hence $\omega(\bar{x})=X$.\par 
	
	\bigskip
	
	Note that $(6)\Longrightarrow (1)$ is obvious.

\end{proof}

\begin{Theorem}	Let $(X,G)$ be a minimal CR-dynamical system. If a closed  $A \subseteq X$ is +invariant in $(X,G)$ then  $A=\emptyset$, or $A =X$.
\end{Theorem}
\begin{proof}
	Suppose that $(X,G)$ is minimal. Let $A$ be a nonempty closed +invariant subset of $X$ and $x\in A$. Clearly, $\cO^+_G(x)\subseteq A$ and as $A$ is closed, it follows that $X = \overline{\cO^+_G(x)}\subseteq A$. 
\end{proof}

The converse part of the above theorem is not generally true.
\begin{Example} We take the same Example as in \cite{mini}.
	
	Let $X=[0,1]$ and $G$ be the union of the following line segments:\par 
	1. the line segment with endpoints $(0,\frac{1}{2})$ and $(1,1)$\par 
	2. the line segment with endpoints $(1,0)$ and $(1,1)$.\par
	
	\bigskip
	
	$(X,G)$ is not minimal as $\mathcal{O}_G^+(0)$ is not dense in $X$. Let $A$ be a nonempty +invariant closed subset of $X$. If $1\in A$, then we are done. As $A\neq \emptyset$, let $x\in A$ with $x<1$. Since $A$ is closed and +invariant so $1\in \overline{\mathcal{O}^+_G(x)} \subseteq A$.
\end{Example}

\begin{Theorem}
	Suppose that $(X,G)$ is a CR-dynamical system. Then the followings are equivalent.
	\begin{enumerate}
		\item $(X,G)$ is minimal.
		\item For every $x\in X,$ there exists $\bar{x}\in T_G^+(x)$ such that $\omega(\bar{x})=X$.
		\item For every $x\in X$, $\omega_G(x)=X$.
	\end{enumerate}
\end{Theorem}
\begin{proof}
	For $(1)\Longrightarrow (2)$ we assume $(X,G)$ be minimal and let $x\in X$. Then from Theorem \ref{minimality}, there exists $\bar{x}\in T_G^+(x)$ such that $\cO_G^{\oplus}(\bar{x})$ is dense in $X$.	 Let $y\in X$ and $\{U_n\}$ be the neighbourhood base at $y$. Then there exists $n_1\in \N$ such that $\pi_{n_1}(\bar{x})\in U_1\setminus \{y\}$.  Take $W_2=U_2\setminus \{\pi_i(\bar{x}):i=0,1,\dots,n_1\}$. Then there exists $n_2 > n_1$ such that $\pi_{n_2}(\bar{x})\in W_2$. Subsequently we get $n_1<n_2<n_3<\dots$ such that $\pi_{n_k}(\bar{x})\in W_k$, where $W_k=U_k\setminus \{\pi_i(\bar{x}):i=0,1,\dots,n_{k-1}\}$ for every $k>1$. Clearly,  $\pi_{n_k}(\bar{x})$ coverges to $y$. Hence $y\in \omega_G(\bar{x})$.\par
	
	Note that $(2)\Longrightarrow (3)\Longrightarrow (1)$ is obivous.

\end{proof}

\begin{Remark} The above theorem gives the answer of problems 5.10 and 5.11 posed in \cite{mini}. Now we can conclude that the concepts of $2^{\oplus}$-minimality, $3^{\oplus}$-minimality, $2^{\omega}$-minimality and $3^{\omega}$-minimality in \cite{mini} are equivalent.\par 

We note that in a point transitive system $(X,G)$ for some $x \in X$, $\omega_G(x)=X$ need not hold. For example if we consider a closed relation $G$ on $[0,1]$ such that $G=\{\frac{1}{2}\}\times [0,1]\cup [0,1]\times \{0\}.$ Any forward orbit of $\frac{1}{2}$ is in the form $\bar{\frac{1}{2}}=(\frac{1}{2},\dots,\frac{1}{2},x,0,0,\dots)$ for some $x\in [0,1]$. Hence $\omega_G(\frac{1}{2})$ contains only $0$ and $\frac{1}{2}$. \end{Remark}

\bigskip

\begin{Definition}
	Let $(X,G)$ be a CR-dynamical system, where $G$ is suitable. Then $(X,G)$ is called \emph{suitably minimal} if for every $x\in X,$ $\cO^{\bullet +}_G(x)$ is dense in $X.$
\end{Definition}

\begin{Remark}
	Note that for a suitable relation $G$, $(X,G)$ is suitably minimal $\iff$ $(X,G, \bullet)$ is minimal.
	
	Also Example \ref{irr} is also suitably minimal.
\end{Remark}

Clearly, suitable minimality implies minimality. But the converse is not true as can be seen in the example below. 

\begin{Example}\label{mini not suit}
 Let $X=[0,1]$ and a closed relation $F$ be defined as $$F=\displaystyle{\bigcup_{n=1}^{\infty}}\big(\big[\frac{1}{n+1},\frac{1}{n}\big]\times \big[\frac{1}{n+1},\frac{1}{n}\big]\big)\bigcup \big\{(0,0),(1,0),(0,1)\big\}.$$ 
	
	Clearly, $(X,G)$ is minimal but with $G$ not being a suitable relation we can not define suitable minimality.
\end{Example}

Note that Example \ref{irr} is suitably minimal but not strongly minimal.

	\begin{Theorem}
		Let $(X,G)$ be a CR-dynamical system, where $G$ is suitable and $g$ is a selection function of $G$ . Then the following conditions are equivalent.\par 
		\begin{enumerate}
			\item $(X,G)$ is suitably minimal.
			\item For every opene  $U\subseteq X$ and for every point $x\in X$, there exists $n\in \mathbb{N}$ such that $G^{\bullet n}(x)\cap U\neq \emptyset$.
			\item For every opene  $U\subseteq X$ and every point $x\in X$, the set $N^{\bullet}(x,U)$ is nonempty.
			\item For every opene set $U\subseteq X$ and for every point $x\in X,$ there exists $n\in \N$ such that there is some sequence $(x_k)$ converges to $x$ with $x_k\in C_g^{\infty}$ and $g^n(x_k)\in U$, for each $k$.
			\end{enumerate}
\end{Theorem}
\begin{proof}
	$(1)\iff (2) \iff (3)$ is similar to Theorem \ref{minimality} and $(2)\iff (3)$ is clear from the fact that for every $n\in \N,$ $G^{\bullet n}=\overline{\{(x,g^n(x)):x\in C_g^{\infty}\}}.$
\end{proof}

\bigskip 

One can naturally ask the following questions in the direction of suitable minimality:\par

\bigskip

\noindent\begin{tabular}{|l|}
	
	\hline\\
	
\textbf{Questions:} \\
1. Is there any example of minimal system $(X,G),$ with $G$ suitable but $(X,G)$\\ is not suitably minimal?\\
2. What is the relation between strong minimality and suitable minimality for \\ a suitable relation?\\
3. How can one define strongly suitably minimal system?\\
4. How can one define invariant set in the suitable dynamics $(X,G,\bullet)?$\\

\hline
\end{tabular}

\bigskip

For dynamics of continuous maps, the concepts of  minimality and almost periodic points are closely related. In CR-dynamical system this presents a more general picture. Here we define two types of almost periodic points and establish a characterisation of minimal sytems with some counter examples.
\begin{Definition}
	Let $(X,G)$ be a CR-dynamical system and $x\in X$. Then\par 
	\begin{enumerate}
		\item $x$ is called a \emph{$G$-orbital almost periodic point} if for every $\varepsilon >0$  there exists $K\in \mathbb{N}$ such that
		
		 $\mathcal{O}_G^{+}(x)\subseteq \displaystyle{\bigcap_{n\in \mathbb{N}}}[B(G^n(x),\varepsilon)\cup B(G^{(n+1)}(x),\varepsilon)\dots \cup B(G^{(n+K)}(x),\varepsilon)]$.\par
		\item $x$ is called a \emph{$G$-trajectorial almost periodic point} if there exists some $\bar{x}\in T_G^+(x)$ such that  for every $\varepsilon >0$ there exists $K\in \mathbb{N}$ with
		
		 $\mathcal{O}_G^{\oplus}(\bar{x})\subseteq \displaystyle{\bigcap_{n\in \mathbb{N}}}[B(\pi_n(\bar{x}),\varepsilon)\cup B(\pi_{(n+1)}(\bar{x}),\varepsilon)\dots \cup B(\pi_{(n+K)}(\bar{x}),\varepsilon)]$. 
		
\end{enumerate}

\end{Definition}

The concepts of $G$-orbital almost periodic points and $G$-trajectorial almost periodic points are mutually independent. We give some counterexamples below:

\begin{Example}
 Let $X=\displaystyle{\bigcup_{n=1}^{\infty}}\big(\{\frac{1}{n}\}\times [0,1]\big) \bigcup \big(\{0\}\times [0,1]\big) \subseteq \R^2$. We define a closed relation 
		
		\bigskip
		
		$G=\big(X\times \{(1,0)\}\big) \bigcup \big\{((\frac{1}{n},0),(\frac{1}{n+1},0)):n\geq 2 \big\} \bigcup \big\{((\frac{1}{n+1},0),(\frac{1}{n},0)):n\geq 3 \big\} \bigcup \big\{((\frac{1}{n},0),(\frac{1}{2},\frac{1}{n})):n\geq 3\big\} \bigcup \big\{((0,0),(\frac{1}{2},0)),((0,0),(0,0)) \big\}.$
		
		\bigskip
		
		 Let $\varepsilon >0$ be given. Then there exists $p\in \N$ such that $\frac{1}{p}<\varepsilon.$ Take $K=2p+3.$  
		 Then $\cO_G^ +((\frac{1}{2},0))\subseteq$ 
		 
		  $ \displaystyle{\bigcap_{n\in \mathbb{N}}}[B(G^n((\frac{1}{2},0)),\varepsilon)\cup B(G^{(n+1)}((\frac{1}{2},0)),\varepsilon)\dots \cup B(G^{(n+K)}((\frac{1}{2},0)),\varepsilon)]$.
		  
		   Hence $(\frac{1}{2},0)$ is a $G$-orbital almost periodic point.
		   
		   \bigskip

		     Let $\bar{x}\in T_G^+((\frac{1}{2},0)).$ If $(1,0)\in \cO^{\oplus}_G(\bar{x})$ that is there exists some $k\in \N$ such that $\pi_{k+1}(\bar{x})=(1,0)$ then
		     $\pi_n(\bar{x})=(1,0)$ for every $n>(k+1).$ Consequently, $(\frac{1}{2},0)\notin B(\pi_n(\bar{x}),\frac{1}{2})$ for every $n>(k+1),$ which implies that $(\frac{1}{2},0)$ is not a $G$-trajectorial almost periodic point. On the other hand if $(1,0)\notin \cO_G^{\oplus}(\bar{x}),$ $\bar{x}$ is of the form $\bar{x}=((\frac{1}{2},0),(\frac{1}{3},0),(\frac{1}{4},0),\dots),$ where $\pi_n(\bar{x})\in \{(\frac{1}{i},0):i\in \N\  \text{and}\  i>2\}$ for all $n>3.$ Hence $(\frac{1}{2},0)\notin B(\pi_n(\bar{x}),\frac{1}{6})$ for every $n>1$.

		    Therefore $(\frac{1}{2},0)$ is not a $G$-trajectorial almost periodic point.
		
	\end{Example}

\begin{Example}  Let $\mathbb{S}$ be an unit circle in $\R^2$ and $\theta: \mathbb{S} \to \mathbb{S}$ be an irrational rotation. Define $f:[1,2]\to [1,2]$ as $f(2-\frac{1}{n})=2-\frac{1}{n+1}$ for every $n\in \N,$ $f(2)=2$ and $f$ is linear elsewhere.

	 \bigskip
	 
	  Let $X=\mathbb{S} \cup \big\{(x,0):x\in [1,2]\big\} \subseteq \R^2$. 
	 We consider a closed relation $G$ on $X$ as
	  $$G=\big\{(x,\theta(x)):x\in \mathbb{S}\big\}\cup \big\{((x,0),(f(x),0)):x\in [1,2]\big\}.$$
	  
	   For $x=(1,0)$ take the forward orbit, $\bar{x}=(x,\theta(x),\theta^2(x),\dots)$. Since $(\mathbb{S},\theta)$ is a minimal system, so for every $\varepsilon>0$ there exists $K\in \N$ such that the set $\{n\in \N:\theta^n(x)\in B(x,\frac{\varepsilon}{2})\}$ is syndetic with gap of length $K.$ Furthermore as $\theta$ is an isometry, we see  that  
	   $$\cO_G^{\oplus}(\bar{x})\subseteq \displaystyle{\bigcap_{n\in \mathbb{N}}}[B(\pi_n(\bar{x}),\varepsilon)\cup B(\pi_{(n+1)}(\bar{x}),\varepsilon)\dots \cup B(\pi_{(n+2K)}(\bar{x}),\varepsilon)].$$
	   
	    Hence $(1,0)$ is a $G$-trajectorial almost periodic point. 
	    
	    On the other hand $(\frac{3}{2},0)\notin G^n((1,0))$ for every $n>2,$ which implies $(\frac{3}{2},0)\notin B(G^n((1,0)),\frac{1}{6})$ for every $n>2.$ This proves that $(1,0)$ is not a $G$-orbital almost periodic point. 
\end{Example}

\begin{Theorem}
	Let $(X,G)$ be a CR-dynamical system.\par $(1)$ If $x$ is a $G$-trajectorial almost periodic point then there exists $\bar{x}\in T_G^+(x)$ such that for every open neighbourhood $U$ of $x$, $N_{\bar{x}}(x,U)$ is syndetic.\par
$(2)$ If $x$ is a $G$-orbital almost periodic point then for every open neighbourhood $U$ of $x$, $N(x,U)$ is syndetic. 
\end{Theorem}
\begin{proof}
	$(1)$ Suppose that $x$ is a $G$-trajectorial almost periodic point and $\bar{x}\in T_G^+(x)$ is the required forward orbit. Let $U$ be an open neighbourhood of $x$ and $\varepsilon >0$ with $B(x,\varepsilon)\subseteq U$. Then there exists $K\in \mathbb{N}$ such that $\mathcal{O}_G^{\oplus}(\bar{x})\subseteq \displaystyle{\bigcap_{n\in \mathbb{N}}}[B(\pi_n(\bar{x}),\varepsilon)\cup B(\pi_{(n+1)}(\bar{x}),\varepsilon)\dots \cup B(\pi_{(n+K)}(\bar{x}),\varepsilon)]$. Consequently, for each $n\in \mathbb{N}$ there exists sone $r_n\in \{0,1,\dots,K\}$ such that $x\in B(\pi_{(n+r_n)},\varepsilon)$. Hence every interval of length $(K+1)$ meets $N_{\bar{x}}(x,U)$ i.e $N_{\bar{x}}(x,U)$ is syndetic.\par 
	$(2)$ The proof is similar to $(1)$.
	
\end{proof}
Converse of the above theorem is not true. We illustrate this fact in the following example. This example also shows that for every open neighbourhood $U$ of $x,$ $N(x,U)$ is syndetic may not imply the existence of a forward orbit $\bar{x}$ with $N_{\bar{x}}(x,U)$ is syndetic for every open neighbourhood $U$ of $x.$\par

\begin{Example}
	Let $X=\displaystyle{\bigcup_{n=1}^{\infty}}\big(\{\frac{1}{n}\}\times [0,1]\big) \bigcup \big(\{0\}\times [0,1]\big) \subseteq \R^2$. We define a closed relation 
	
	\bigskip
	
	$G=\big(X\times \{(1,0)\}\big) \bigcup \big\{((\frac{1}{n},0),(\frac{1}{n+1},0)):n\geq 2 \big\} \bigcup \big\{((\frac{1}{n},0),(\frac{1}{2},\frac{1}{n})):n\geq 3\big\} \bigcup \big\{((0,0),(\frac{1}{2},0)),((0,0),(0,0)) \big\}.$
	
	\bigskip

	 One can observe that for every open neighbourhood $U$ of $(\frac{1}{2},0),$ $N((\frac{1}{2},0),U)$ is syndetic. But $(\frac{1}{2},0)$ is not a $G$-orbital almost periodic point as if we take $\varepsilon=\frac{1}{12},$ then $(\frac{1}{3},0)\notin B(G^k((\frac{1}{2},0)),\epsilon)$ for every $k\geq 2.$

	  Furthermore, let $\bar{x}\in T_G^+((\frac{1}{2},0)).$ If $(1,0)\in \cO^{\oplus}_G(\bar{x}),$ there exists some $k\in \N$ such that $\pi_n(\bar{x})=(1,0)$ for every $n>k.$ Consequently, $(\frac{1}{2},0)\notin B(\pi_n(\bar{x}),\frac{1}{2})$ for every $n>k,$ which implies that $N_{\bar{x}}((\frac{1}{2},0),B((\frac{1}{2},0),\frac{1}{2}))$ is not syndetic. On the other hand if $(1,0)\notin \cO_G^{\oplus}(\bar{x}),$ $\bar{x}$ is of the form $\bar{x}=((\frac{1}{2},0),(\frac{1}{3},0),(\frac{1}{4},0),\dots),$ where $\pi_n(\bar{x})\in \{(\frac{1}{i},0):i\in \N\  \text{and}\  i>2\}$ for all $n>3.$ Hence $(\frac{1}{2},0)\notin B(\pi_n(\bar{x}),\frac{1}{6})$ for every $n>1.$ Hence the set $N_{\bar{x}}((\frac{1}{2},0),B((\frac{1}{2},0),\frac{1}{6}))$ is not syndetic. 
	 
	 Consequently, $(\frac{1}{2},0)$ is not $G$-trajectorial almost periodic point.
\end{Example}

\bigskip
\noindent\begin{tabular}{|l|}
	\hline\\
	
\textbf{Question.} Is it true that if $N_{\bar{x}}(x,U)$ is syndetic for every open  neighbourhood  $U$ \\ of $x$ for some forward orbit $\bar{x}$ then $x$ is $G$-trajectorial almost periodic point? \\

\hline
\end{tabular}

\bigskip

\begin{Theorem}\label{AP}
	Suppose that $(X,G)$ is a CR-dynamical system and $x\in X$ is a $G$-trajectorial almost periodic point. Then there exists a forward orbit $\bar{x}\in T_G^+(x)$ such that  $\overline{\mathcal{O}_G^{\oplus}(\bar{x})}$ is minimal.
\end{Theorem}
\begin{proof}
	Suppose that $x$ is a $G$-trajectorial almost periodic point and $\bar{x}\in T_G^+(x)$ is the corresponding forward orbit. First of all we will show that for every $y\in \overline{\mathcal{O}_G^{\oplus}(\bar{x})}$ there exists a forward orbit $\bar{y}\in T_G^+(y)$ such that $\mathcal{O}_G^{\oplus}(\bar{y})\subseteq \overline{\mathcal{O}_G^{\oplus}(\bar{x})}$. If $y\in \mathcal{O}_G^{\oplus}(\bar{x})$ then $y=\pi_t(\bar{x})$ for some $t\in \mathbb{N}$ and so we can get a forward orbit $\bar{y}\in T_G^+(y)$ such that $\pi_n(\bar{y})=\pi_{(t+n-1)}(\bar{x})$ for every $n\in \mathbb{N}$. If $y\in \overline{\mathcal{O}_G^{\oplus}(\bar{x})}\setminus \mathcal{O}_G^{\oplus}(\bar{x})$, then there exists a sequence $(\pi_{t_i}(\bar{x}))\to y$. Since $X$ is compact, one can observe that for every $n\in \mathbb{N}$, $(\pi_{(t_i+n)}(\bar{x}))_{i\in \mathbb{N}}$ has a convergent subsequence that converges to some $y_n\in \overline{\mathcal{O}_G^{\oplus}(\bar{x})}$. Without any loss of generality passing to a subsequence if necessary, we can assume that $(\pi_{(t_i+n)}(\bar{x}))_{i\in \mathbb{N}}$ converges to $y_n$. From closedness of $G^n$ we can say that $y_n\in G^n(y)$. We can now construct a forward orbit $\bar{y}=(y=y_0,y_1,y_2,\dots)\in T_G^+(y)$ where $\pi_n(\bar{y})=y_{(n-1)}$ and so $\mathcal{O}_G^{\oplus}(\bar{y})\subseteq \overline{\mathcal{O}_G^{\oplus}(\bar{x})}$.  Now we will show that $\mathcal{O}_G^{\oplus}(\bar{y})$ is dense in $\overline{\mathcal{O}_G^{\oplus}(\bar{x})}$. Let $\varepsilon >0$ be given and $n\in \mathbb{N}$. Then there exists $N\in \mathbb{N}$ such that $\mathcal{O}_G^{\oplus}(\bar{x})\subseteq \displaystyle{\bigcap_{r\in \mathbb{N}}}[B(\pi_r(\bar{x}),\frac{\varepsilon}{2})\cup B(\pi_{(r+1)}(\bar{x}),\frac{\varepsilon}{2})\dots \cup B(\pi_{(r+N)}(\bar{x}),\frac{\varepsilon}{2})]$. If $y=\pi_t(\bar{x})$ for some $t\in \N$, then there exists $k\in\{0,1,\dots,N\}$ with $\pi_n(\bar{x})\in B(\pi_{(t+k)}(\bar{x}),\frac{\varepsilon}{2})$ and so $\pi_k(\bar{y})\in B(\pi_n(\bar{x}),\frac{\varepsilon}{2})$. If $y\in \overline{\mathcal{O}_G^{\oplus}(\bar{x})}\setminus \mathcal{O}_G^{\oplus}(\bar{x})$, that is $(\pi_{t_i}(\bar{x}))\to y$ then for each $i\in \mathbb{N}$, there exists $k_i\in \{0,1,\dots,N\}$ such that $\pi_n(\bar{x})\in B(\pi_{(t_i+k_i)}(\bar{x}),\frac{\varepsilon}{2})$. Without any loss of generality passing to a subsequence we can assume that $\pi_n(\bar{x})\in B(\pi_{(t_i+k)}(\bar{x}),\frac{\varepsilon}{2})$ for some $k\in \{0,1,\dots,N\}$. We have already observed that $(\pi_{(t_i+k)}(\bar{x}))$ converges to $y_k$. So $y_k\in B(\pi_n(\bar{x}),\varepsilon)$. Consequently, $B(\pi_n(\bar{x}),\varepsilon)\cap \mathcal{O}_G^{\oplus}(\bar{y})$. Since $\varepsilon >0$ and $n\in \mathbb{N}$ are choosen arbitrarily, $\mathcal{O}_G^{\oplus}(\bar{y})$ is dense in $\overline{\mathcal{O}_G^{\oplus}(\bar{x})}$.
\end{proof}
\begin{Example}
	Let $X=\displaystyle{\bigcup_{n=1}^{\infty}}\big(\{\frac{1}{n}\}\times [0,1]\big) \bigcup \big(\{0\}\times [0,1]\big) \subseteq \R^2$. We define a closed relation 
	
	\bigskip
	
  $G=\big\{((\frac{1}{n},x),(\frac{1}{n+1},x)):n\in \N\ \text{and}\  x\in [0,1]\big\}\bigcup \big\{((0,x),(0,x)):x\in [0,1]\big\}\bigcup \big\{((\frac{1}{2},0),(1,0))\big\}.$ 
	\bigskip

		Here $(1,0)$ is $G$-trajectorial almost periodic and has a forward orbit
		
		 $((1,0),(\frac{1}{2},0),(1,0),\dots),$ being periodic it is a minimal set.
		
		On the other hand $(1,0)$ is also a $G$-orbital almost periodic point. Indeed let $\varepsilon >0,$ be given. Then there exists $k\in \N,$ such that $(\frac{1}{m},0)\in B((\frac{1}{n},0),\varepsilon)$ for every $m,n\geq k.$ Take $K=k+1.$ 	
		Then
		
		\noindent $\mathcal{O}_G^{+}((1,0))\subseteq \displaystyle{\bigcap_{n\in \mathbb{N}}}[B(G^n((1,0)),\varepsilon)\cup B(G^{(n+1)}((1,0)),\varepsilon)\dots \cup B(G^{(n+K)}((1,0)),\varepsilon)]$. 
		
		But $\overline{\cO_G^+((1,0))}$ is not minimal. Since $\cO_G^+((0,0))$ is not dense in $\overline{\cO_G^+((1,0))}.$
\end{Example}

\begin{Theorem}
Let $(X,G)$ be a CR-dynamical system. Suppose that there exist a point $x$ and a forward orbit $\bar{x}$ such that $N_{\bar{x}}(x,U)$ is syndetic for every open neighbourhood $U$ of $x.$ Then for every $y\in \overline{\cO_G^{\oplus}(\bar{x})},$ $x\in \omega_G(y).$
\end{Theorem}
\begin{proof}
	Let $y\in \overline{\cO_G^{\oplus}(\bar{x})}.$ Then proceeding as the same way as in Theorem \ref{AP}, we can construct a forward orbit $\bar{y}=(y=y_0,y_1,y_2,\dots)$ such that $\cO_G^{\oplus}(\bar{y})\subseteq \overline{\cO_G^{\oplus}(\bar{x})}.$ Let $y=\pi_t(\bar{x})$ for some $t\in \N$. For every open set $U$ containing $x$ the set $N_{\bar{x}}(x,U)$ is syndetic and so infinite, which implies that $x\in \omega(\bar{y}).$\par 
	
	Suppose that $y\in \overline{\mathcal{O}_G^{\oplus}(\bar{x})}\setminus \mathcal{O}_G^{\oplus}(\bar{x})$, and $(\pi_{t_i}(\bar{x}))\to y.$ Let $\varepsilon>0$ be given. Then there exists $N\in \N$ such that for each $i\in \mathbb{N}$, there exists $k_i\in \{0,1,\dots,N\}$ with $\pi_{(t_i+k_i)}(\bar{x})\in B(x,\frac{\varepsilon}{2})$. Without any loss of generality passing to a subsequence we can assume that for each $i,$ $\pi_{(t_i+k)}(\bar{x})\in B(x,\frac{\varepsilon}{2})$ for some $k\in \{0,1,\dots,N\}$. From the proof of Theorem \ref{AP}, $ (\pi_{(t_i+k)}(\bar{x}))$ converges to $y_k$. So $y_k\in B(x,\varepsilon)$. Consequently, $x\in \omega(\bar{y}).$ 
\end{proof}

\begin{Theorem}
	Minimality and strong minimality of closed relation is a residual property.
\end{Theorem}
\begin{proof} We prove the same in three steps.
	\begin{enumerate}
		\item \underline{For factors:} 
		
		Suppose that $h:(X,G)\to (Y,F)$ is a factor map and $(X,G)$ is minimal. Let $y\in Y$, then there exists a $x\in X$  with $h(x)=y$. Let $U$ be opene in $Y.$ Since $(X,G)$ is minimal, there exists $n\in \N,$ such that $G^n(x)\cap h^{-1}(U)\neq \emptyset,$ which implies that $h(G^n(x))\cap U \neq \emptyset$ and so $F^n(h(x))\cap U\neq \emptyset.$ Hence $\displaystyle{\bigcup_{k=1}^{\infty}}F^k(y)$ is dense in $Y$ and thus $(Y,F)$ is minimal.\par 
	
	Suppose that $(X,G)$ is strongly minimal and $h:(X,G)\to (Y,F)$ is a factor map. Let $y\in Y$, then there exists a $x\in X$  with $h(x)=y$. Let $\bar{y}=(y= y_1,y_2,\dots)\in T_F^+(y)$ be a forward orbit. Now $y_2\in F(h(x)))$ and so $ y_2\in h(G(x)),$ which implies that there exists $x_2\in G(x)$ with $y_2=h(x_2).$ In this process we can construct a forward orbit $(x=x_1,x_2,\dots)\in T_G^+(x)$ such that for each $k,$ $h(x_k)=y_k$. Let $U$ be  opene in $Y.$ Since $(X,G)$ is strongly minimal, there exists $n\in \N,$ such that $x_n \in h^{-1}(U)\Rightarrow h(x_n)\in U.$ Hence $\cO_F^{\oplus}(\bar{y})$ is dense in $Y$ and so $(Y,F)$ is strongly minimal. 
	\vspace{0.2cm}
		\item \underline{For irreducible extensions:}

		Suppose that $h:(X,G)\to (Y,F)$ is an irreducible epimorphism and $(Y,F)$ is minimal. Let $x\in X$ be a point and $U$ be opene in $X.$ Then there exists an opene $V \subseteq Y$ such that $h^{-1}(V)\subseteq U.$ Since $(Y,F)$ is minimal, there exists $n\in \N,$ such that $F^n(h(x))\cap V\neq \emptyset,$ which implies that $h(G^n(x))\cap V \neq \emptyset$ and then $G^n(x)\cap h^{-1}(V)\neq \emptyset.$ Hence $\displaystyle{\bigcup_{k=1}^{\infty}}G^k(x)$ is dense in $X$ and so $(X,G)$ is minimal.\par 
		
		Suppose that $(Y,F)$ is strongly minimal. Let $x\in X$ be a point and $\bar{x}=(x=x_1,x_2,\dots)\in T_G^+(x)$ be a forward orbit. Now for each $k>1,$ $h(x_k)\in h(G(x_{k-1})) \Rightarrow h(x_k)\in F(h(x_{k-1})).$ Hence we can construct a forward orbit $(y_1,y_2,\dots)$ such that for each $k,$ $h(x_k)=y_k$. Let us consider $U$ be an opene in $X.$ Then there exists an opene $V$ in $Y$ such that $h^{-1}(V)\subseteq U.$ Since $(Y,F)$ is strongly minimal, there exists $n\in \N,$ such that $h(x_n)\in V \Rightarrow x_n\in U.$ Hence $\cO_G^{\oplus}(\bar{x})$ is dense in $X$ and so $(X,G)$ is strongly minimal.
		\vspace{0.2cm}
		\item  \underline{For inverse limits:}
		
		Let $\{(X_n,F_n):n\in \mathbb{N}\}$ be a family of CR-dynamical system together with epimorphisms $\{h_i:(X_{i+1},F_{i+1})\to (X_i,F_i):i\in \mathbb{N}\}$ and $(X,F)$ be the associated inverse limit and $p_k:(X,F)\to (X_k,F_k)$ be the kth coordinate projection. Assume that each $(X_i,F_i)$ is minimal. Let $x\in X$ and $U$ be opene in $X.$ Then from Lemma \ref{inverselimit}, there exist some $k$ and opene $V$ in $X_k$ such that $p_k^{-1}(V)\subseteq U.$ Since $(X_k,F_k)$ is minimal, there exists $n\in \N$ such that $F_k^n(p_k(x))\cap V \neq \emptyset,$ which implies that $p_k(F^n(x))\cap V\neq \emptyset.$ Thus $F^n(x)\cap U\neq \emptyset.$ Hence $(X,F)$ is minimal.\par
		
		Let us assume that each $(X_i,F_i)$ is strongly minimal. Let $\mathbf{x}=(x_k)\in X$ and consider a forward orbit $\bar{\mathbf{x}}=(\mathbf{x}=\mathbf{x}_1,\mathbf{x}_2,\dots,\mathbf{x}_i,\dots)\in T_F^+(\mathbf{x}),$ where for each $i\in \N,$ $\mathbf{x}_i=(x^i_k).$ It is evident that for each $n\in \N,$ $\bar{x_n}=(x_n=x^1_n,x^2_n,\dots)\in T_{F_n}^+(x_n)$ and $\cO_{F_n}^{\oplus}(\bar{x_n})$ is dense in $X_n.$ Let $U$ be opene in $X.$ Then from Lemma \ref{inverselimit}, there exists $k\in \N$ and opene $U_k\in X_k$ such that $p_k^{-1}(U_k)\subseteq U.$ Furthermore, there exists some $r\in \N$ with $x^r_k\in U_k \implies p_k(\mathbf{x}_r)\in U_k,$ which implies that $\mathbf{x}_r\in U.$ Therefore $(X,F)$ is strongly minimal. 
	\end{enumerate}

\end{proof}

We see that for a CR-dynamical system $(X,G)$ we have:

\bigskip
\noindent\begin{tabular}{|l|}
	\hline\\
	
	$\begin{array}{c}
		\text{Strongly Minimal}\\
		\text{Suitably Minimal}\\
	\end{array}$  $ \Longrightarrow$ Minimal $\Longrightarrow$ Strongly Point Transitive\\
	
	\hline
\end{tabular}

\bigskip

\bigskip

\subsection{Strongly Transitive, Very Strongly Transitive and   suitable forms}

\begin{Definition} Let $(X,G)$ be a CR-dynamical system.
	
	$(X,G)$ is called \emph{Strongly Transitive}  if for every opene set $U\subseteq X$, 
	$$\displaystyle{\bigcup_{n=1}^{\infty}}G^n(U)=X.$$ 
	
	$(X,G)$ is called \emph{Very Strongly Transitive}  if for every opene set $U\subseteq X$, there is $N\in \mathbb{N}$ such that $$\displaystyle{\bigcup_{n=1}^{N}}G^n(U)=X.$$

\end{Definition}

Note that a very strongly transitive $(X,G)$ is always strongly transitive.

\bigskip

In the dynamics of continuous map, a minimal system is always very strongly transitive. But in this may not be true in general for CR-dynamical systems. We consider a variation of Example \ref{mini not suit},

\begin{Example}
	 Let $X=[0,1]$ and a closed relation $F$ be defined as $$F=\big (X\times \{1\}\big)\displaystyle{\bigcup_{n=1}^{\infty}}\big(\big[\frac{1}{n+1},\frac{1}{n}\big]\times \big[\frac{1}{n+1},\frac{1}{n}\big]\big)\bigcup \big\{(0,0), (1,0)\big\}.$$

	Clearly, $(X,G)$ is minimal and strongly transitive. Moreover, for every point $x\in  [\frac{1}{k+1},\frac{1}{k}],$ one can easily see that $G^{-(k+1)}(x)=X$ and $G^{-2}(0)=X$. But it is not very strongly transitive as for every $n\in \mathbb{N}$, $\displaystyle{\bigcup_{r=1}^n}G^r((\frac{1}{2},1))\neq X$. 
\end{Example}

Consider a variation of Example \ref{everything},

\begin{Example} \label{vst}
	Let $X=[0,1]$ and a closed relation 
	$$G=\big(\displaystyle{\bigcup_{n=1}^{\infty}}[\frac{1}{2^{n}},\frac{1}{2^{n-1}}]\times [\frac{1}{2^{n}},\frac{1}{2^{n-1}}] \big) \cup \big(\displaystyle{\bigcup_{n=0}^{\infty}}[\frac{1}{3.2^{n}},\frac{1}{3.2^{n-1}}]\times [\frac{1}{3.2^{n}},\frac{1}{3.2^{n-1}}] \big) \cup \big\{ \{\frac{1}{2}\} \times [0,1] \big\} \cup\{(0,0)\}.$$ 
	
	Then $(X,G)$ is very strongly transitive but not minimal since $\cO_G^+ ((0,0)) = \{(0,0)\}$.
	
\end{Example}

\begin{Theorem} \label{st}
	Suppose that $(X,G)$ is a CR-dynamical system. Then the following conditions are equivalent.\par 
	$(1)$ $(X,G)$ is strongly transitive.\par 
	$(2)$ For every opene set $U\subseteq X$ and for every point $x\in X$, there exists $n\in \mathbb{N}$ such that $x\in G^{n}(U)$.\par 
	$(3)$ For every opene set $U\subseteq X$ and every point $x\in X$, the set $N(U,x)$ is nonempty.\par 
	$(4)$ For every opene set $U\subseteq X$ and every point $x\in X$, the set $N(U,x)$ is infinite.\par
	$(5)$ For every $x\in X$, $\mathcal{O}_G^{-}(x)$ is dense in $X$. (This condition is called as $3^{\ominus}$-minimality in \cite{mini}.)\par 
	$(6)$ For every $x\in X$, $G^{-1}(x)\neq \emptyset$ and $\displaystyle{\bigcup_{n=1}^{\infty}}G^{-n}(x)$ is dense in $X$.\par 
	$(7)$ If $E$ is nonempty and -invariant, then $E$ is dense in $X$.
	
	\bigskip
	
	If $(X,G)$ is strongly  transitive, then $(X,G)$ is topologically transitive.
\end{Theorem}
\begin{proof}
	The proof of $(1)\iff (2)\iff (3)$ is clear from definition and $(4)\Longrightarrow
	(3)$ is obvious.\par 
	$(1) \Longrightarrow (5):$ Let $(X,G)$ be strongly transitive and $x\in X$. Clearly, $G(X)=X$. Now for every opene set $U\subseteq X$, we have $\displaystyle{\bigcup_{n=1}^{\infty}}G^{n}(U)=X$. Then there is some $y\in U$ and $k\in \mathbb{N}$ such that $x\in G^{k}(y)$. Then we can construct a backward orbit $\bar{x}=(x,\dots,y,\dots)\in \displaystyle{\bigstar_{i=1}^{\infty}}G^{-1}$ such that $\pi_{k+1}(\bar{x})=y$, which implies that $\mathcal{O}_G^{-}(x)\cap U \neq \emptyset$. Hence $\mathcal{O}_G^{-}(x)$ is dense in $X$.\par 
	$(5)\Rightarrow (1):$ On the contrary suppose that $(X,G)$	is not strongly transitive. Then there is an opene set $U\subseteq X$ such that $\displaystyle{\bigcup_{n=1}^{\infty}}G^n(U)\neq X$. Let $x\notin \displaystyle{\bigcup_{n=1}^{\infty}}G^n(U)$. Since $\mathcal{O}_G^{-}(x)$ is dense in $X$, there is a backward orbit $\bar{x}^{-1}\in T_G^{-}(x)$ and $N\in \mathbb{N}$ such that $\pi_{(N+1)}(\bar{x})\in U\setminus \{x\}$. Then we get $(x=x_1,\dots,x_{(N+1)})\in \displaystyle{\bigstar_{i=1}^{N}}G^{-1}$ with $x_{(N+1)}\in U$, which implies that $x\in G^N(U)$.\par 
	$(3)\Longrightarrow (4)$ On the contrary suppose that there are $x\in X$ and an opene set $U\subseteq X$ such that $N(U,x)$ is finite. Let the maximum element of $N(U.x)$ is $N$. Then there exists $y\in U$ such that $x\in G^N(y)$. Let $k\in N(U,y)$, which implies that there is $z\in U$ such that $y\in G^k(z)$ and so $x\in G^{N+k}(z)$. Hence $N+k\in N(U,x)$, which is a contradiction.\par 
	$(5)\Longrightarrow (6)$. It is clear that $G(X)=X$. The another part is clear from that fact that for every $n\in \mathbb{N}$ and $\bar{x}^{-1}\in T_G^{-}(x)$, $\pi_{n+1}(\bar{x})\in G^{-n}(x)$.\par 
	$(6)\Longrightarrow (5)$ The proof is clear from the fact that $y\in G^{-n}(x)$ implies that there is a backward orbit $\bar{x}^{-1}\in T_G^{-}(x)$ such that $y = \pi_{(n+1)}(\bar{x})$.\par 
	$(5)\iff (7)$. $\mathcal{O}_G^{-}(x)$ is -invariant and if $x\in E$ and $E$ is -invariant then $\mathcal{O}_G^{-}(x)\subseteq E$.
\end{proof}

\begin{Theorem} For a very strongly transitive CR-dynamical system  $(X,G)$, the set $N(U,x)$ is syndetic for every  opene set $U \subseteq X$ and every point $x \in X$. 
\end{Theorem}

\begin{proof}
	If $X = \bigcup \limits_{n = 1}^{N} \ G^n(U)$ then for every
	$k \in \N, \ X = G^k(X) = \bigcup \limits_{n = k+1}^{N+k} \ G^n(U)$.
	Thus, for every $x \in X$, the set $N(U,x)$ meets every interval of length $N$ in $\N$. \end{proof}

\begin{Corollary} If $(X,G)$ is very strongly transitive then for any opene $U, V \subseteq X$,	the set $N(U,V)$ is syndetic. \end{Corollary}

\begin{proof} If $x \in V$, then $N(U,x) \subseteq N(U,V)$.	
\end{proof} 

\bigskip

\begin{Definition}
	Let $(X,G)$ be a CR-dynamical system where $G$ is a suitable relation. Then  
	
-	$(X,G)$ is called \emph{suitably strongly transitive}   if for every  opene  $U \subseteq X$,  $\bigcup \limits_{n = 1}^\infty G^{\bullet n}(U) = X$.
	
-	$(X,G)$ is called \emph{suitably very strongly transitive} if for every  opene  $U \subseteq X$, there exists  $N \in \N$ such that $\bigcup \limits_{n = 1}^N G^{\bullet n}(U) = X$.
	
\end{Definition}

Clearly suitably very strongly transitive implies suitably strongly transitive. We note that Example \ref{irr} is suitably very strongly transitive. 

\bigskip

\begin{tabular}{|l|}
	\hline\\
	
	\textbf{Question:} Does suitably minimal $\Longrightarrow$ suitably very strongly transitive?\\
	
	\hline
\end{tabular}

\bigskip

\begin{Theorem}
	Suppose that $(X,G)$ is a CR-dynamical system where $G$ is suitable. Then the following conditions are equivalent.\par 
	$(1)$ $(X,G)$ is suitably strongly transitive.\par 
	$(2)$ For every opene set $U\subseteq X$ and for every point $x\in X$, there exists $n\in \mathbb{N}$ such that $x\in G^{\bullet n}(U)$.\par 
	$(3)$ For every opene set $U\subseteq X$ and every point $x\in X$, the set $N^\bullet(U,x)$ is nonempty.\par 
	\par
	 
	\bigskip
	
	If $(X,G)$ is suitably strongly  transitive, then $(X,G)$ is suitably topologically transitive.
\end{Theorem}

The proof is similar to the proof of Theorem \ref{st} and we skip refurbishing the details in this proof.

\bigskip

Note that for CR-dynamical system $(X,G)$ we have:

\bigskip
\noindent\begin{tabular}{|l|}
	\hline\\
	
{{\small
		$\begin{array}{ccccc}
			\text{suitably very strongly trnasitive}  & \Longrightarrow & \text{ suitably strongly transitive} & \Longrightarrow &  \text{ suitably transitive}\\
			\Downarrow &    & \Downarrow &   & \Downarrow\\
			\text{ very strongly transitive}  & \Longrightarrow & \text{  strongly transitive} & \Longrightarrow &  \text{  transitive}
		\end{array}$}}	\\
	
	\hline
\end{tabular}

\bigskip

\bigskip

\subsection{Exact Transitive, Strongly Exact Transitive and  suitable forms} 	

\begin{Definition} A CR-dynamical system $(X,G)$ is called \emph{exact} if for every pair of opene subsets $U,V \subseteq X$
	there exists $n \in \N$ such  that $G^n(U) \cap G^n(V) \ \neq \emptyset$.
	
	The system $(X,G)$ is \emph{fully exact} if for every pair of opene subsets $U,V \subseteq X$
	there exists $n \in \N$ such  that $(G^n(U) \cap G^n(V))^{\circ} \ \neq \emptyset$.
\end{Definition}

\begin{Definition} Let $(X,G)$ be a CR-dynamical system.
	
-	$(X,G)$ is called \emph{Exact Transitive}  if for every pair of opene sets $U,V\subseteq X$, $\displaystyle{\bigcup_{n=1}^{\infty}} (G^n(U)\cap G^n(V))$ is dense in $X$.\par 
	
-	$(X,G)$ is called \emph{Strongly Exact Transitive}  if for every pair of opene sets $U,V \subseteq X$, $\displaystyle{\bigcup_{n=1}^{\infty}} (G^n(U)\cap G^n(V)) =X$.
\end{Definition}

As observed in \cite{vt-1} these concepts are far from transitivity. But these give rise to alternate forms of transitivity.

\begin{Theorem} \label{set} Let $(X,G)$ be a CR-dynamical system. Then	
	the following conditions are equivalent:
	\begin{enumerate}
		\item  The system is strongly exact transitive.
		
		\item For every pair of opene sets $U, V \subset X$, the set \ $\bigcup_{n \in \N} \ (G \times G)^n(U \times V)$
		contains the diagonal $\Delta$.
		
		\item For every $x \in X$,  $\cO_{G \times G}^-(x,x)$ is dense in $X \times X$.
	\end{enumerate}
	
	If $(X,G)$ is strongly exact transitive then it is exact transitive.
\end{Theorem}

\begin{proof}All three conditions say that for every $x \in X$ and opene $U, V \subseteq X$ there exists $n \in \N$ such that
	$x \in G^n(U)\cap G^n(V)$.
	
\end{proof}

\begin{Remark}
	Observe that Example \ref{everything} is exact transitive but not strongly exact transitive whereas Example \ref{vst} is  strongly exact transitive.
\end{Remark}

\begin{Theorem} Let $(X,G)$ be a CR-dynamical system. Then,
	
	(a)  If $(X,G)$ is exact transitive then it is topologically transitive and  exact.
	
	(b) If $(X,G)$ is topologically transitive and fully exact then it is exact transitive.
	
	(c) If $(X,G)$ is strongly exact transitive then it is fully exact.
	
	(d) If $(X,G)$ is strongly exact transitive then it is strongly transitive.
\end{Theorem}

\begin{proof} (a) Follows immediately from definition.
	
	(b)  Assume $(X,G)$ is topologically transitive and fully exact.
	For  opene  $U,V, W \subseteq X$  there exists   $k \in \N$ such that $ (G^k(U) \cap G^k(V))^\circ \neq \emptyset$. Transitivity gives  some $m \in \N$ such that $G^m((G^k(U) \cap G^k(V))^\circ) \cap W \neq \emptyset$, i.e. there exists $x \in W$ with $x \in G^{m+k}(U) \cap G^{m+k}(V)$.
	Thus $\bigcup \limits_{n=1}^\infty  (G^n(U) \cap G^n(V))$  is dense in $ X $.
	
	(c)  This is an easy consequence of the Baire Category Theorem.
	
	(d) Follows from definition.
	
\end{proof}

\bigskip

\begin{Remark}
	This shows that exact transitivity is a slight strengthening of the conjunction of exactness and topological transitivity for closed relations similar to continuous functions. 
\end{Remark}

\bigskip 

\begin{Definition} Consider the    CR-dynamical system $(X,G, \bullet)$.
	
-	$(X,G)$ is called \emph{Suitably Exact Transitive}  if for every pair of opene sets $U,V\subseteq X$, $\displaystyle{\bigcup_{n=1}^{\infty}} (G^{\bullet n}(U)\cap G^{\bullet n}(V))$ is dense in $X$.\par 
	
-	$(X,G)$ is called \emph{Suitably Strongly Exact Transitive}  if for every pair of opene sets $U,V \subseteq X$, $\displaystyle{\bigcup_{n=1}^{\infty}} (G^{\bullet n}(U)\cap G^{\bullet n}(V)) =X$.
\end{Definition}

\begin{Theorem} \label{et} For suitable $G$, let $(X,G)$ be a suitably exact transitive CR-dynamical system. Then it is suitably weakly mixing.\end{Theorem}
\begin{proof}
	For opene $U, V  \subseteq X$, we see that there exists  $n \in \N$ such that $U \cap G^{\bullet n}(U) \neq \emptyset$ and $ G^{\bullet n}(U)\cap G^{\bullet n}(V) \neq \emptyset$. And so by Theorem \ref{wm}, $(X,G)$ is suitably weakly mixing.
\end{proof}

\bigskip

Note that for CR-dynamical system $(X,G)$ we have:

\bigskip
\noindent\begin{tabular}{|l|}
	\hline\\
	
	{{\small
			$\begin{array}{ccccc}
				\text{suitably  strongly exact transitive}  & \Longrightarrow & \text{ suitably exact transitive} & \Longrightarrow &  \text{ suitably weakly mixing}\\
				\Downarrow &    & \Downarrow &   & \Downarrow\\
				\text{ strongly exact transitive}  & \Longrightarrow & \text{  exact transitive} &  &  \text{  weakly mixing}\\
				\Downarrow &    & \Downarrow &   & \Downarrow\\
				\text{ strongly  transitive}  & \Longrightarrow & \text{  transitive} &  &  \text{  transitive}
			\end{array}$}}	\\
	
	\hline
\end{tabular}

\bigskip

\bigskip 

\bigskip

\subsection{Strongly and Suitably Strongly Product Transitive} Let $(X,G)$ be a CR-dynamical system.

\begin{Definition}
	$(X,G)$ is called \emph{Strongly Product Transitive}  if for every positive integer $k$, the product system $(X^k,G^{[k]})$ is strongly transitive.
\end{Definition}

\begin{Remark}
	If $(X,G)$ is strongly product transitive then it is strong exact transitive.
\end{Remark}

\bigskip 

\begin{Definition}
	$(X,G)$ is called \emph{Suitably Strongly Product Transitive}  if for every positive integer $k$, the product system $(X^k,G^{[k]}, \bullet)$ is strongly transitive.
\end{Definition}

\begin{Theorem} For a CR-dynamical system $(X,G)$, with $G$ suitable,  the
	following are equivalent:
	\begin{enumerate}
		\item  The system is suitably strongly product transitive.

		\item The collection of subsets  $\{ N^\bullet(U,x) : x \in X $ and $U$ opene in $X \}$ of $\N$ has the finite intersection
		property (or equivalently it generates a filter of subsets of $\N$).
	\end{enumerate}
	
	If $(X,G)$ is suitably strongly product transitive then it is suitably strong exact transitive.
\end{Theorem}

We have for $(X,G)$:

\bigskip
\noindent\begin{tabular}{|l|}
	\hline\\
	
	{{\small
			$\begin{array}{ccccc}
				\text{suitably  strongly product transitive}  & \Longrightarrow & \text{ strongly product transitive}   \\
				\Downarrow &    & \Downarrow    \\
				\text{ suitably strongly exact transitive}  & \Longrightarrow &  \text{ strongly exact transitive}   \\
			\end{array}$}}	\\
	
	\hline
\end{tabular}

\bigskip

\bigskip

\subsection{Locally Eventually Onto and its suitable version} Let $(X,G)$ be a CR-dynamical system.

\begin{Definition}
	$(X,G)$ is called Locally Eventually Onto  if for every opene $U\subseteq X$, there exists $N\in \mathbb{N}$, $G^N(U)=X$, and so $G^n(U)=X$ for all $n\geq N$.
\end{Definition}

\begin{Theorem}
	
	If $(X,G)$ is locally eventually onto then it is strongly product transitive and topologically mixing.
\end{Theorem}

We consider a modification of Example \ref{everything},

\begin{Example} 
	Let $X=[0,1]$ and a closed relation 
	$$G=\big(\displaystyle{\bigcup_{n=1}^{\infty}}[\frac{1}{2^{n}},\frac{1}{2^{n-1}}]\times [\frac{1}{2^{n}},\frac{1}{2^{n-1}}] \big) \cup \big(\displaystyle{\bigcup_{n=0}^{\infty}}[\frac{1}{3.2^{n}},\frac{1}{3.2^{n-1}}]\times [\frac{1}{3.2^{n}},\frac{1}{3.2^{n-1}}] \big)$$
	\centerline{$ \cup \big\{\{\frac{1}{n}, n \in \N\} \times [0,1]\big\}\cup \big\{ \{0\} \times [0,1] \big\}.$}
	
	Then $(X,G)$ is locally eventually onto.
	
\end{Example}

\begin{Remark}
	It is clear that if $(X,G)$ is locally eventually onto then $N(U,x)$ is co-finite in $\N$ for all opene
	$U \subseteq X$ and all $x \in X$.  Converse need not be true.
\end{Remark} 

\begin{Definition}
	For a suitable $G$, the CR-dynamical system $(X,G)$ is called \emph{Suitably Locally Eventually Onto}  if for every opene $U\subseteq X$, there exists $N\in \mathbb{N}$, $G^{\bullet N}(U)=X$, and so $G^{\bullet n}(U)=X$ for all $n\geq N$.
\end{Definition}

\begin{Theorem}
	
	If $(X,G)$ is suitably locally eventually onto then it is suitably strongly product transitive and suitably topologically mixing.
\end{Theorem}

\begin{Remark}
	It is clear that if $(X,G)$ is suitably locally eventually onto then $N^{\bullet}(U,x)$ is co-finite in $\N$ for all opene
	$U \subseteq X$ and all $x \in X$.  Converse need not be true.
\end{Remark} 

\bigskip
\noindent\begin{tabular}{|l|}
	\hline\\
	
	{{\small
			$\begin{array}{ccccc}
				\text{suitably  locally eventually onto}  &  \begin{array}{l}
						\Longrightarrow \text{Suitably Mixing} \\
						\Longrightarrow \text{Suitably Strongly Product Transitive}\\
						\end{array}\\
				\Downarrow &    \\
				\text{ locally eventually onto}  &  \begin{array}{l}
						\Longrightarrow \text{Mixing} \\
						\Longrightarrow \text{Strongly Product Transitive}\\
						\end{array}
			\end{array}$}}	\\
	
	\hline
\end{tabular}

\bigskip

\bigskip

We find that various implications hold for these different kinds of transitivities. We come accross several questions for which we have no answer as of now. It is still to be settled if the reverse implications in all the above cases do not hold.

\bigskip

\noindent\textbf{Acknowledgement:} The first named author is thankful to UGC (University Grants Commission of
India, (Ref: 1127/(CSIR-UGC NET DEC. 2017))), for giving Junior Research Fellowship during the tenure of which this work has been done.

\end{document}